\documentclass[11pt]{amsart}

\usepackage{amssymb,amsthm,amsmath}
\usepackage[numbers,sort&compress]{natbib}
\usepackage{color}
\usepackage{graphicx}
\usepackage{tikz}
\usepackage[colorlinks,
            linkcolor=blue,
            anchorcolor=blue,
            citecolor=blue
            ]{hyperref}

\let\newpf\proof \let\proof\relax 
\newenvironment{pf}{\newpf[\proofname]}{\qed\endtrivlist}

\newcommand{\ba}{\overline{A}}

\def\be{\begin{equation}}
\def\ee{\end{equation}}

\def\ba{{\begin{align}}}
\def\ea{{\end{align}}}

\def\bm{\begin{matrix}}
\def\em{\end{matrix}}

\def\0{{\mathbf 0}}

\newtheorem{Theorem}{Theorem}[section]
\newtheorem*{Theorem*}{Theorem 1.1}
\newtheorem{Lemma}{Lemma}[section]
\newtheorem{Proposition}{Proposition}[section]

\newtheorem{Remark}{Remark}[section]

\newtheorem{Claim}{Claim}

\numberwithin{equation}{section}

\theoremstyle{definition}

\renewcommand{\mod}{\operatorname{mod}}

\newcommand{\C}{{\mathbb C}}

\newcommand{\N}{{\mathbb N}}

\newcommand{\R}{{\mathbb R}}
\newcommand{\T}{{\mathbb T}}

\newcommand{\Z}{{\mathbb Z}}

\newcommand{\la}{\langle}
\newcommand{\ra}{\rangle}

\def\B0{{\bold{0}}}

\catcode`\@=12

\def\Empty{}
\newcommand\oplabel[1]{
  \def\OpArg{#1} \ifx \OpArg\Empty {} \else
    \label{#1}
  \fi}

\newcommand{\comm}[1]{}
\newcommand{\comment}[1]{}

\begin{document}

\title[Absolutely Continuous Spectrum]{The absolutely continuous spectrum of finitely differentiable quasi-periodic Schr\"{o}dinger operators}

\author{Ao Cai} \address{
Chern Institute of Mathematics and LPMC, Nankai University, Tianjin 300071, China; and Departmento de Matem\'{a}tica and CMAFCIO, Faculdade de Ci\^{e}ncias, Universidade de Lisboa, Portugal.
}

\email{acai@fc.ul.pt; godcaiao@126.com}

%\author{Qi Zhou} \address{
%}
%
%\email{}
\setcounter{tocdepth}{1}

\begin{abstract}
We prove that the quasi-periodic Schr\"{o}dinger operator with a finitely differentiable potential has purely absolutely continuous spectrum for all phases if the frequency is Diophantine and the potential is sufficiently small in the corresponding $C^k$ topology. This is based on a refined quantitative $C^{k,k_0}$ almost reducibility theorem which only requires a quite low initial regularity ``$k>14\tau+2$'' and much of the regularity ``$k_0\leq k-2\tau-2$'' is conserved in the end, where $\tau$ is the Diophantine constant of the frequency.
\end{abstract}

\maketitle
\section{Introduction}
In this paper, we shall consider the quasi-periodic Schr\"odinger operator $H_{V,\alpha,\theta}$ with a finitely differentiable potential:
\begin{equation}\label{1.1}
(H_{V,\alpha,\theta}x)_n=x_{n+1}+x_{n-1}+V(\theta+n\alpha)x_{n}, n\in \Z,
\end{equation}
where $\theta\in \T^d=\R^d/(2\pi\Z)^d$ is called the phase, $\alpha \in \R^d$ is called the frequency satisfying $\la m,\alpha\ra \notin 2\pi \Z$ for any $m\in \Z^d$ different from zero, and $V\in C^k(\T^d,\R)$ is called the potential, $k,d\in \N^+$. A typical example is the almost Mathieu operator $H_{2\lambda \cos,\alpha,\theta}$:
$$
(H_{2\lambda \cos,\alpha,\theta}x)_n=x_{n+1}+x_{n-1}+2\lambda \cos(\theta+n\alpha)x_{n}, n\in \Z,
$$
where $\lambda$ is called the coupling constant.

Quasi-periodic Schr\"odinger operators come from solid-state physics, showing the influence of an external magnetic field on the electrons of a crystal \cite{Da,MJ}. It is related to quasi-crystal which is between conductor and insulator. As is very well known, the absolutely continuous spectrum corresponds to conductor. That is, the absolutely continuous spectrum of a quasi-periodic Schr\"odinger operator is the set of energies at which the described physical system exhibits transport. Moreover, the absolutely continuous spectrum is the part of spectrum that has the best stability properties under small perturbation and its existence has strong implications, including an Oracle Theorem that predicts the potential, as shown by Remling \cite{rem}. It is thus an important question to ask whether the quasi-periodic Schr\"odinger operator has (even purely) absolutely continuous spectrum.

Recall that $\alpha \in\R^d$ is called {\it Diophantine} if there are $\kappa>0$ and $\tau>d$ such that $\alpha \in {\rm DC}(\kappa,\tau)$, where
\begin{equation}\label{dio}
{\rm DC}(\kappa,\tau):=\left\{\alpha \in\R^d:  \inf_{j \in \Z}\left|\la n,\alpha  \ra - 2\pi j \right|
> \frac{\kappa}{|n|^{\tau}},\quad \forall \, n\in\Z^d\backslash\{0\} \right\}.
\end{equation}
Here we denote
$$
\lvert n\rvert=\lvert n_1\rvert+\lvert n_2\rvert+ \cdots + \lvert n_d\rvert,
$$
and
$$
\la n,\alpha \ra =n_1\alpha_1+n_2\alpha_2+\cdots+n_d\alpha_d.
$$
Denote ${\rm DC}=\bigcup_{\kappa,\tau} {\rm DC}(\kappa,\tau)$, which is of full Lebesgue measure.

Our main theorem is the following:
\begin{Theorem}\label{main}
Assume $\alpha\in {\rm DC}(\kappa, \tau)$, $V\in C^k(\T^d, \R)$ with $k>35\tau+2$. If $\lambda$ is sufficiently small, then $H_{\lambda V,\alpha,\theta}$ has purely absolutely continuous spectrum for all $\theta$.
\end{Theorem}

\begin{Remark}
We do not claim the optimality of the lower bound of ``$k$'' but we point out that some kind of regularity is essential for the existence of $($purely$)$ absolutely continuous spectrum, as will be stated later. For some technical reason, we require ``$k$'' to be larger than $35\tau+2$ instead of $14\tau+2$.
\end{Remark}
It appears that the existence of (purely) absolutely continuous spectrum depends sensitively on the arithmetic properties of the frequency. Recently, Avila and Jitomirskaya \cite{AJ3} constructed super-Liouvillean $\alpha\in\T^2$ such that for typical analytic potential, the corresponding quasi-periodic Schr\"{o}dinger operator has no absolutely continuous spectrum. Relatively, Hou-Wang-Zhou \cite{HWZ} showed that there exists super-Liouvillean $\alpha\in\T^2$ such that for small analytic potential, the corresponding quasi-periodic Schr\"{o}dinger operator has absolutely continuous spectrum. Moreover, they proved that if $\alpha\in\T^d$ with $\alpha$ being weak-Liouvillean and the potential is small enough, the absolutely continuous spectrum exists.

When $d=1$, things can be characterized much more explicitly thanks to Avila's fantastic global theory of analytic Schr\"{o}dinger operators \cite{avila}. He showed that typical one-frequency operators have only point spectrum in the supercritical region, and absolutely continuous spectrum in the subcritical region. It seems that the effect of $\alpha$'s arithmetic properties is weaker in one-frequency case, but we point out that the proofs of absolutely continuous spectrum are quite different according to the arithmetic assumptions on $\alpha$. Focusing on the almost Mathieu operator $H_{2\lambda \cos,\alpha,\theta}$, there is a famous conjecture: Simon \cite{simon2} (Problem $6$) asked whether AMO has purely absolutely continuous spectrum for all $0<\lvert \lambda\rvert<1$, all phases and all frequencies. This conjecture was first proved for Diophantine $\alpha$ and almost every $\theta$ by Jitomirskaya \cite{J}, whose approach follows Aubry duality and localization theory. About ten years later, two key advances happened. On one hand, Avila and Jitomirskaya \cite{AJ2} established so-called quantitative duality to prove Simon's conjecture for Diophantine $\alpha$ and all $\theta$. On the other hand, Avila and Damanik \cite{AD2} used periodic approximation and Kotani theory to prove the conjecture for Liouvillean $\alpha$ and almost every $\theta$. The complete solution to Simon's problem was given by Avila \cite{A01}. He distinguished the whole proof into two parts: when $\beta=0$ (the subexponential regime, see \cite{A01}), the proof relied on almost reducibility results developed in \cite{AJ2}; when $\beta>0$ (the exponential regime), he improved the periodic approximation method developed in \cite{AD2}.

Another important factor which influences the existence of absolutely continuous spectrum is the regularity of the potential. In the analytic topology, Dinaburg and Sinai \cite{DS} proved that $H_{V,\alpha,\theta}$ has absolutely continuous spectrum component for all $\theta$ in the perturbative regime ($V$ being analytically small and the smallness depends on $\alpha$) by reducibility theory. Later, Eliasson \cite{Eli92} improved the KAM scheme and showed that $H_{V,\alpha,\theta}$ has purely absolutely continuous spectrum for all $\theta$ in the same setting. By ``non-perturbative reduction to perturbative regime'', Eliasson's result was extended to the non-perturbative regime by Avila-Jitomirskaya \cite{AJ2} and Hou-You \cite{houyou}.

However, when it comes to the finitely differentiable topology, there are few results on this issue. In this sense, our theorem is constructive and it mainly extends Eliasson's \cite{Eli92} result to the finitely differentiable case. We emphasize that assuming some kind of continuity or higher regularity of the potential is necessary, not only for ``purely'' absolutely continuous spectrum, but also for the existence of an absolutely continuous spectrum component. In $C^0$ topology, Avila and Damanik \cite{AD} proved that for one-dimensional Schr\"{o}dinger operators with ergodic continuous potentials, there exists a generic set of such potentials such that the corresponding operators have no absolutely continuous spectrum. Moreover, by Gordon's Lemma, Boshernitzan and Damanik \cite{BD} proved that for generic ergodic continuous potentials, the corresponding operators have purely singular continuous spectrum.

Now, let us show the main strategy of our paper in the following. As is known to all, Schr\"{o}dinger operator $H_{V,\alpha,\theta}$ is closely related to Schr\"{o}dinger cocycle $(\alpha,A)$ where
$$
A(\theta)=S_E^{V}(\theta)=\begin{pmatrix} E- V(\theta) & -1 \\ 1 & 0 \end{pmatrix},
$$
since the solution of $H_{V,\alpha,\theta}x=Ex$ satisfies
$$
A(\theta+n\alpha)\begin{pmatrix} x_{n} \\ x_{n-1} \end{pmatrix}=\begin{pmatrix} x_{n+1} \\ x_{n} \end{pmatrix}.
$$
Therefore, we can use reducibility method to analyze the dynamics of $C^k$ quasi-periodic linear cocycle $(\alpha, A)\in \T^d \times C^k(\T^d,SL(2,\R))$ and then study the spectral properties of the corresponding operator. This approach, which was first developed in \cite{Eli92},  has been proved to be very fruitful \cite{A01,AJ2,AYZ1,AYZ,LYZZ}. Readers are invited to consult You's 2018 ICM survey \cite{you} for more related achievements. In this paper, we acquire a finer quantitative $C^k$ almost reducible theorem by distinguishing resonances from non-resonances more precisely. For simplicity, let us introduce its qualitative version on the Schr\"{o}dinger cocycle (for the quantitative one which works for general $C^k$ $SL(2,\R)$-valued cocycles, see Theorem \ref{thm3}).

\begin{Theorem}\label{thm1.2}
Let $\alpha \in {\rm DC}(\kappa,\tau)$, $V\in C^{k}(\T^{d},\R)$ with $k>14\tau+2$. If $\lambda$ is sufficiently small, then for any $E\in\R$, $(\alpha,S_E^{\lambda V})$ is $C^{k,k_0}$ almost reducible with $k_0\leq k-2\tau-2$.
\end{Theorem}
\begin{Remark}
If we change the assumption into $k>17\tau+2$, we can further prove the $\frac{1}{2}$-H\"{o}lder continuity of the Lyapunov exponent and the integrated density of states $($see Theorem \ref{thm1} and Theorem \ref{cor1}$)$, which greatly reduces the initial regularity requirement of $k\geqslant 550\tau$ in \cite{CCYZ}. Moreover, compared with \cite{CCYZ}, we obtain a quite good upper bound of the remainder $k_0\leq k-2\tau-2$ instead of $k_0\leq k/6$ for $C^{k,k_0}$ almost reducibility.
\end{Remark}

Finally, we point out that our main idea of proving the purely absolutely continuous spectrum follows that of Avila \cite{A01} in the subexponential regime. There are two important aspects. One is that we need to obtain a modified quantitative $C^k$ almost reducibility theorem which was originally established by Cai-Chavaudret-You-Zhou \cite{CCYZ}. It will provide us with fine estimates on the conjugation map, the constant part and the perturbation in each KAM step. The other is, we need to stratify the spectrum of $H_{V,\alpha,\theta}$ by the rotation number of $(\alpha,A)$. Once they are done, we will be able to have a good control of the growth of the transfer matrix (see (\ref{transfer}) for definition) on each hierarchical spectrum part. The proofs left are standard by the theorems of Gilbert-Pearson \cite{GP} and Avila \cite{A01}.

\section{Preliminaries}
For a bounded
analytic (possibly matrix valued) function $F(\theta)$ defined on $\mathcal{S}_h:= \{ \theta=(\theta_1,\dots, \theta_d)\in \mathbb{C}^d\ |\  \forall 1\leqslant i\leqslant d, \ | \Im \theta_i |< h \}$, let
$
|F|_h=  \sup_{\theta\in \mathcal{S}_h } \| F(\theta)\| $ and denote by $C^\omega_{h}(\T^d,*)$ the
set of all these $*$-valued functions ($*$ will usually denote $\R$, $sl(2,\R)$,
$SL(2,\R)$). We denote $C^\omega(\T^d,*)=\cup_{h>0}C^\omega_{h}(\T^d,*)$, and set $C^{k}(\T^{d},*)$ to be the space of $k$ times differentiable with continuous $k$-th derivatives functions. The norm is defined as
$$
\lVert F \rVert _{k}=\sup_{\substack{
                             \lvert k^{'}\rvert\leqslant k,
                             \theta \in \T^{d}
                          }}\lVert \partial^{k^{'}}F(\theta) \rVert.
$$

\subsection{Conjugation and reducibility}
Given two cocycles $(\alpha,A_1)$, $(\alpha,A_2)\in \T^d    \times C^{\ast}(\T^d,SL(2,\R))$, ``$\ast$'' stands for ``$\omega$'' or ``$k$'',  one says that they are $C^{\ast}$ conjugated if there exists $Z\in C^{\ast}(2\T^d, SL(2,\R))$, such that $$
Z(\theta+\alpha)A_1(\theta)Z^{-1}(\theta)=A_2(\theta).
$$
Note that we need to define $Z$ on the $2\T^d=\R^d/(4\pi\Z)^d$ in order to make it still real-valued.

An analytic cocycle $(\alpha, A)\in \T^d    \times C^{\omega}_h(\T^d,SL(2,\R))$ is said to be almost reducible if there exist a sequence of conjugations $Z_j\in C^{\omega}_{h_j}(2\T^d, SL(2,\R))$, a sequence of constant matrices $A_j\in SL(2,\R)$ and a sequence of small perturbation $f_j \in C^{\omega}_{h_j}(\T^d, sl(2,\R))$ such that
$$
Z_j(\theta+\alpha)A(\theta)Z_j(\theta)^{-1}=A_j e^{f_j(\theta)}
$$
with
$$
\lvert f_j(\theta)\rvert_{h_j}\rightarrow 0, \ \ j\rightarrow \infty.
$$
Furthermore, we call it weak $(C^{\omega})$ almost reducible if $h_j \rightarrow 0$ and we call it strong $(C^{\omega}_{h,h'})$ almost reducible if $h_j\rightarrow h'>0$. We say $(\alpha, A)$ is $C^{\omega}_{h,h'}$ reducible if there exist a conjugation map $\tilde{Z}\in C^{\omega}_{h'}(2\T^d, $ $SL(2,\R))$ and a constant matrix $\tilde{A} \in SL(2,\R)$ such that
$$
\tilde{Z}(\theta+\alpha)A(\theta)\tilde{Z}(\theta)^{-1}=\tilde{A}(\theta).
$$

In order to avoid repetition, we give an equivalent definition of $C^k$ (almost) reducibility in the following.

A finitely differentiable cocycle $(\alpha,A)$ is said to be $C^{k,k_1}$ almost reducible, if $A\in C^k(\T^d,SL(2,\R))$ and the $C^{k_1}$-closure of its $C^{k_1}$ conjugacies contains a constant. Moreover, we say $(\alpha,A)$ is  $C^{k,k_1}$ reducible, if $A\in C^k(\T^d,SL(2,\R))$ and its $C^{k_1}$ conjugacies contain a constant.

\subsection{Analytic approximation}
Assume $f \in C^{k}(\T^{d},sl(2,\R))$. By Zehnder \cite{zehnder}, there exists a sequence $\{f_{j}\}_{j\geqslant 1}$, $f_{j}\in C_{\frac{1}{j}}^{\omega}(\T^{d},sl(2,\R))$ and a universal constant $C^{'}$, such that
\begin{eqnarray} \nonumber \lVert f_{j}-f \rVert_{k} &\rightarrow& 0 , \quad  j \rightarrow +\infty, \\
\label{2.1}\lvert f_{j}\rvert_{\frac{1}{j}} &\leqslant& C^{'}\lVert f \rVert_{k}, \\   \nonumber \lvert f_{j+1}-f_{j} \rvert_{\frac{1}{j+1}} &\leqslant& C^{'}(\frac{1}{j})^k\lVert f \rVert_{k}.
\end{eqnarray}
Moreover, if $k\leqslant \tilde{k}$ and $f\in C^{\tilde{k}}$, then properties $(\ref{2.1})$ hold with $\tilde{k}$ instead of $k$. That means this sequence is obtained from $f$ regardless of its regularity (since $f_{j}$ is the convolution of $f$ with a map which does not depend on $k$).

\subsection{Rotation number and degree}
Assume that $A\in C^0(\T^d,SL(2,\R))$ is homotopic to identity. It introduces the projective skew-product $F_A:\T^d \times \mathbb{S}^1 \rightarrow \T^d \times \mathbb{S}^1$ with
$$
F_A(x,\omega):=\big(  x+\alpha, \frac{A(x)\cdot \omega}{\lvert A(x)\cdot \omega\rvert}\big),
$$
which is also homotopic to identity. Thus we can lift $F_A$ to a map $\tilde{F}_A:\T^d\times \R\rightarrow \T^d\times \R$ of the form $\tilde{F}_A(x,y)=(x+\alpha,y+\psi(x,y))$, where for every $x \in \T^d$, $\psi(x,y)$ is $2\pi\Z$-periodic in $y$. The map $\psi:\T^d \times \R \rightarrow \R$ is called a lift of $A$. Let $\mu$ be any probability measure on $\T^d\times \R$ which is invariant by $\tilde{F}_A$, and whose projection on the first coordinate is given by Lebesgue measure. The number
\begin{equation}\label{rho}
\rho_{(\alpha,A)}:=\frac{1}{(2\pi)^d}\int_{\T^d\times \R} \psi(x,y)d\mu(x,y)\mbox{mod}\,2\pi\Z
\end{equation}
does not depend on the choices of the lift $\psi$ or the measure $\mu$. It is called the {\it fibered rotation number} of cocycle $(\alpha,A)$ (readers can consult \cite{JM} for more details).

Let
$$
R_{\phi}:=\begin{pmatrix} \cos\phi & -\sin\phi \\ \sin\phi & \cos\phi \end{pmatrix},
$$
if $A\in C^0(\T^d,SL(2,\R))$ is homotopic to $\theta\rightarrow R_{\la n,\theta\ra}$ for some $n\in \Z^d$, then we call $n$ the {\it degree} of $A$ and denote it by deg$A$. Moreover, \begin{equation}\label{degree}\deg(AB)=\deg A+\deg B.\end{equation}

Note that the fibered rotation number is invariant under real conjugacies which are homotopic to identity. More generally, if the cocycle $(\alpha,A_1)$ is conjugated to $(\alpha,A_2)$ by $B\in C^0(2\T^d, SL(2,\R))$, i.e. $B(\cdot +\alpha)A_1(\cdot)B^{-1}(\cdot)=A_2(\cdot)$, then
\begin{equation}\label{degpro}
\rho_{(\alpha,A_2)}=\rho_{(\alpha,A_1)}+\frac{\la \deg B,\alpha\ra}{2}.
\end{equation}
\section{Dynamical estimates: almost reducibility}\label{sec3}
In this section, we establish the modified quantitative $C^{k,k_0}$ almost reducibility for finitely differentiable quasi-periodic $SL(2,\R)$ cocycles, which will be applied to control the growth of corresponding Schr\"{o}dinger cocycles.

Consider the $C^k$ quasi-periodic $SL(2,\R)$ cocycle:
$$
(\alpha,Ae^{f(\theta)}):\T^{d}\times\R^{2} \rightarrow \T^{d}\times\R^{2};(\theta,v)\mapsto (\theta+\alpha,Ae^{f(\theta)}\cdot v),
$$
where
$A\in SL(2,\R), \, f\in C^k(\T^{d},sl(2,\R)), \, d\in \N^+$ and $\alpha\in {\rm DC}(\kappa,\tau)$. We will first analyze the analytic approximating cocycles $\{(\alpha,Ae^{f_{j}(\theta)})\}_{j\geqslant 1}$ and then obtain the estimates of $(\alpha,Ae^{f(\theta)})$ by analytic approximation \cite{zehnder}.

We would like to mention that we will carry out our proof in the framework of $SU(1,1)$ and $su(1,1)$ to make it more explicit. Recall that $sl(2,\R)$ is isomorphic to $su(1,1)$, which consists of matrices of the form
$$
\begin{pmatrix} it & v\\ \bar{v} & -it \end{pmatrix}
$$
with $t\in \R$, $v\in\C$.
The isomorphism between them is given by $A\rightarrow MAM^{-1}$, where
$$
M=\frac{1}{1+i}\begin{pmatrix} 1 & -i\\ 1 & i \end{pmatrix}
$$
and a simple calculation yields
$$
M\begin{pmatrix} x & y+z\\ y-z & -x \end{pmatrix}M^{-1}=\begin{pmatrix} iz & x-iy\\ x+iy & -iz \end{pmatrix},
$$
where $x,y,z\in \R$. $SU(1,1)$ is the corresponding Lie group of $su(1,1)$.
\subsection{Notations}
In the following subsections, parameters $\rho,\epsilon,N,\sigma$ will be fixed; one will refer to the situation where there exists $n_\ast$ with $0<\lvert n_\ast\rvert \leqslant N$ such that
$$
\inf_{j \in \Z}\lvert 2\rho- \la n_\ast,\alpha\ra - 2\pi j\rvert< \epsilon^{\sigma},
$$
as the ``resonant case'' (for simplicity, we just write ``$\lvert 2\rho - \la n_\ast,\alpha\ra \rvert$'' to represent the left side, same for $\lvert \la n_\ast,\alpha\ra \rvert$). The integer vector $n_\ast$ will be referred to as a ``resonant site''.
Resonances are linked to a useful decomposition of the space $\mathcal{B}_r:=C^{\omega}_{r}(\T^{d},su(1,1))$.

Assume that for given $\eta>0$, $\alpha\in \R^{d}$ and $A\in SU(1,1)$, we have a decomposition $\mathcal{B}_r=\mathcal{B}_r^{nre}(\eta) \bigoplus\mathcal{B}_r^{re}(\eta)$ satisfying that for any $Y\in\mathcal{B}_r^{nre}(\eta)$,
\begin{equation}\label{space}
A^{-1}Y(\theta+\alpha)A\in\mathcal{B}_r^{nre}(\eta),\,\lvert A^{-1}Y(\theta+\alpha)A-Y(\theta)\rvert_r\geqslant\eta\lvert Y(\theta)\rvert_r.
\end{equation}
And let $\mathbb{P}_{nre}$, $\mathbb{P}_{re}$ denote the standard projections from $\mathcal{B}_r$ onto $\mathcal{B}_r^{nre}(\eta)$ and $\mathcal{B}_r^{re}(\eta)$ respectively.

Then we have the following crucial lemma which helps us remove all the non-resonant terms:

\begin{Lemma}\label{lem2}\cite{CCYZ}\cite{houyou}
Assume that $A\in SU(1,1)$, $\epsilon\leqslant (4\lVert A\rVert)^{-4}$ and  $\eta \geqslant 13\lVert A\rVert^2{\epsilon}^{\frac{1}{2}}$. For any $g\in \mathcal{B}_r$ with $|g|_r \leqslant \epsilon$,  there exist $Y\in \mathcal{B}_r$ and $g^{re}\in \mathcal{B}_r^{re}(\eta)$ such that
$$
e^{Y(\theta+\alpha)}(Ae^{g(\theta)})e^{-Y(\theta)}=Ae^{g^{re}(\theta)},
$$
with $\lvert Y \rvert_r\leqslant \epsilon^{\frac{1}{2}}$ and $\lvert g^{re}\rvert_r\leqslant 2\epsilon$.
\end{Lemma}

\begin{Remark}\label{rem2}
In the inequality ``$\eta \geqslant 13\lVert A\rVert^2{\epsilon}^{\frac{1}{2}}$'', ``$\frac{1}{2}$'' is sharp due to the quantitative Implicit Function Theorem \cite{BerBia,deimling}. The proof only relies on the fact that $\mathcal{B}_r$ is a Banach space, thus it also applies to $C^k$ and $C^0$ topology. One can refer to the appendix of \cite{CCYZ} for details.
\end{Remark}

\subsection{Analytic KAM Theorem} In this part, we focus on the analytic quasi-periodic $SL(2,\R)$ cocycle:
$$
(\alpha,Ae^{f(\theta)}):\T^{d}\times\R^{2} \rightarrow \T^{d}\times\R^{2};(\theta,v)\mapsto (\theta+\alpha,Ae^{f(\theta)}\cdot v),
$$
where
$A\in SL(2,\R), \, f\in C^{\omega}_{r}(\T^{d},sl(2,\R))$ with $r>0, d\in \Z^{+}$,
and $\alpha\in {\rm DC}(\kappa,\tau)$.
Note that $A$ has eigenvalues $\{e^{i\rho},e^{-i\rho}\}$ with $\rho \in \R\cup i\R$. We formulate our quantitative analytic KAM theorem as follows.

\begin{Theorem}\label{prop1}
Let $\alpha\in {\rm DC}(\kappa,\tau)$, $\kappa,r>0$, $\tau>d$, $\sigma<\frac{1}{6}$.
Suppose that $A\in SL(2,\R)$, $f\in C^{\omega}_{r}(\T^{d},sl(2,\R))$.  Then for any $r'\in (0,r)$, there exist constants $c=c(\kappa,\tau,d)$, $D> \frac{2}{\sigma}$ and $\tilde{D}=\tilde{D}(\sigma)$ such that if
\begin{equation}\label{estf}
\lvert f \rvert_r\leqslant\epsilon \leqslant \frac{c}{\lVert A\rVert^{\tilde{D}}}(r-r')^{D\tau},
\end{equation}
then there exist $B\in C^{\omega}_{r'}(2\T^{d},SL(2,\R))$, $A_{+}\in SL(2,\R)$ and $f_{+}\in C^{\omega}_{r'}(\T^{d},$
$sl(2,\R))$ such that
$$
B(\theta+\alpha)(Ae^{f(\theta)})B^{-1}(\theta)=A_{+}e^{f_+(\theta)}.
$$
More precisely, let $N=\frac{2}{r-r'} \lvert \ln \epsilon \rvert$, then we can distinguish two cases:
\begin{itemize}
\item $($Non-resonant case$)$   if for any $n\in \Z^{d}$ with $0<\lvert n \rvert \leqslant N$, we have
$$
\lvert 2\rho - \la n,\alpha\ra \rvert\geqslant \epsilon^{\sigma},
$$
then
$$\lvert B-Id\rvert_{r'}\leqslant \epsilon^{\frac{1}{2}} ,\   \ \lvert f_{+}\rvert_{r'}\leqslant \epsilon^{3-\sigma}.$$
and
$$\lVert A_+-A\rVert\leqslant 2\lVert A\rVert\epsilon.$$
\item $($Resonant case$)$ if there exists $n_\ast$ with $0<\lvert n_\ast\rvert \leqslant N$ such that
$$
\lvert 2\rho- \la n_\ast,\alpha\ra \rvert< \epsilon^{\sigma},
$$
then
\begin{flalign*}
\lvert B \rvert_{r'} &\leqslant 8(\frac{\lVert A\rVert}{\kappa})^{\frac{1}{2}}(\frac{2}{r-r'} \lvert \ln \epsilon \rvert)^{\frac{\tau}{2}}\times\epsilon^{\frac{-r'}{r-r'}},\\ \lVert B\rVert_0 &\leqslant 8(\frac{\lVert A\rVert}{\kappa})^{\frac{1}{2}}(\frac{2}{r-r'} \lvert \ln \epsilon \rvert)^{\frac{\tau}{2}},\\
 \lvert f_{+}\rvert_{r'}&\leqslant \frac{2^{5+\tau}\lVert A\rVert\lvert \ln\epsilon\rvert^{\tau}}{\kappa(r-r')^{\tau}}\epsilon e^{-N'(r-r')}(N')^de^{Nr'}\ll \epsilon^{100}, \, N'> 2N^2.
\end{flalign*}
Moreover, $A_+=e^{A''}$ with $\lVert A''\rVert \leqslant 2\epsilon^{\sigma}$, $A''\in sl(2,\R)$. More accurately, we have
$$
MA''M^{-1}=\begin{pmatrix} it & v\\ \bar{v} & -it \end{pmatrix}
$$
with $\lvert t\rvert \leqslant \epsilon^{\sigma}$ and
$$
\lvert v \rvert\leqslant \frac{2^{4+\tau}\lVert A\rVert\lvert \ln\epsilon\rvert^{\tau}}{\kappa(r-r')^{\tau}}\epsilon e^{-\lvert n_{\ast}\rvert r}.
$$
\end{itemize}
\end{Theorem}

\begin{pf}
For the readers who are quite familiar with the analytic KAM scheme, this proof can be skipped since the structure is similar to that in \cite{CCYZ}. But the estimates here are sharp compared with those in \cite{CCYZ} (see Remark \ref{rem3.3}), so we prefer to provide the detailed proof for self-containedness. Let us prove this theorem in $SL(2,\R)$'s isomorphic group: $SU(1,1)$.
We distinguish the proof into two cases:

\textbf{Non-resonant case:}
For $0<\lvert n \rvert \leqslant N=\frac{2}{r-r'} \lvert \ln \epsilon \rvert$, we have
\begin{equation}\label{est1}
\lvert 2\rho - \la n,\alpha\ra \rvert\geqslant \epsilon^{\sigma};
\end{equation}
by $(\ref{estf})$ with $D>\frac{2}{\sigma}$, we have
\begin{equation}\label{est2}
\left \lvert \la n,\alpha\ra \right \rvert \geqslant \frac{\kappa}{\left \lvert n \right \rvert ^{\tau}}\geqslant \frac{\kappa}{\left \lvert N \right \rvert ^{\tau}}\geqslant \epsilon^{\frac{\sigma}{2}}\geqslant \epsilon^{\sigma}.
\end{equation}

It is well known that $(\ref{est1})$ and $(\ref{est2})$ are the conditions which are used to overcome the small denominator problem in KAM theory.

Define
\begin{equation}\label{lambdaN}
\Lambda_N=\{f\in C^{\omega}_{r}(\T^{d},su(1,1))\mid f(\theta)=\sum_{k\in \Z^{d},0<\lvert k \rvert<N}\hat{f}(k)e^{i\la k,\theta\ra}\}.
\end{equation}

Our goal is to solve the cohomological equation
$$
Y(\theta+\alpha)A-AY(\theta)=A(-\mathcal{T}_Nf(\theta)+\hat{f}(0)),
$$
i.e.
\begin{equation}\label{coho}
A^{-1}Y(\theta+\alpha)A-Y(\theta)=-\mathcal{T}_Nf(\theta)+\hat{f}(0).
\end{equation}
Take the Fourier transform for $(\ref{coho})$ and compare the corresponding Fourier coefficients of the two sides. By $(\ref{est1})$ (apply it twice to solve the off-diagonal) along with $(\ref{est2})$ (apply it once to solve the diagonal), we obtain that if $Y\in\Lambda_N$, then
$$
\lvert Y(\theta)\rvert_r \leqslant \epsilon^{-3\sigma}\lvert \mathcal{T}_Nf(\theta)-\hat{f}(0)\rvert_r,
$$
which gives
\begin{equation}\label{crucial}
\lvert A^{-1}Y(\theta+\alpha)A-Y(\theta)\rvert_r\geqslant\epsilon^{3\sigma}\lvert Y(\theta)\rvert_r.
\end{equation}
Moreover, we have $A^{-1}Y(\theta+\alpha)A \in \Lambda_N$ by $(\ref{lambdaN})$. For $\eta=\epsilon^{3\sigma}$, we define $\mathcal{B}_r^{nre}(\epsilon^{3\sigma})$ by $(\ref{space})$, then we have $\Lambda_N \subset \mathcal{B}_r^{nre}(\epsilon^{3\sigma})$.

Since $\epsilon^{3\sigma}\geqslant 13\lVert A\rVert^2\epsilon^{\frac{1}{2}}$ (it holds by $\sigma$ being smaller than $\frac{1}{6}$ and $\tilde{D}$ depending on $\sigma$), by Lemma \ref{lem2} we have $Y\in \mathcal{B}_r$ and $f^{re}\in \mathcal{B}_r^{re}(\epsilon^{3\sigma})$ such that
$$
e^{Y(\theta+\alpha)}(Ae^{f(\theta)})e^{-Y(\theta)}=Ae^{f^{re}(\theta)},
$$
with $\lvert Y \rvert_r\leqslant \epsilon^{\frac{1}{2}}$ and
\begin{equation}\label{fre}
\lvert f^{re}\rvert_r\leqslant 2\epsilon.
\end{equation}
By $(\ref{lambdaN})$
$$
(\mathcal{T}_N{f^{re}})(\theta)=\hat{f}^{re}(0), \ \ \lVert \hat{f}^{re}(0)\rVert \leqslant 2\epsilon,
$$
and
\begin{flalign}\label{ep}
\lvert (\mathcal{R}_N{f^{re}})(\theta)\rvert_{r'}&= \lvert \sum_{\lvert n \rvert>N}\hat{f}^{re}(n)e^{i\la n,\theta\ra}\rvert_{r'}\\
\notag &\leqslant 2\epsilon e^{-N(r-r')}(N)^{d}\\
\notag &\leqslant 2\epsilon\cdot\epsilon^{2}\cdot \frac{1}{4}\epsilon^{-\sigma}\\
\notag &=\frac{1}{2}\epsilon^{3-\sigma}.
\end{flalign}
Moreover, we can compute that
$$
e^{\hat{f}^{re}(0)+\mathcal{R}_N{f^{re}}(\theta)}=e^{\hat{f}^{re}(0)}(Id+e^{-\hat{f}^{re}(0)}\mathcal{O}(\mathcal{R}_N{f^{re}}))=e^{\hat{f}^{re}(0)}e^{f_+(\theta)},
$$
by $(\ref{ep})$, we have
$$
\lvert f_+(\theta)\rvert_{r'}\leqslant 2\lvert \mathcal{R}_N{f^{re}(\theta)}\rvert_{r'} \leqslant \epsilon^{3-\sigma}.
$$
Finally, if we denote
$$
A_+=Ae^{\hat{f}^{re}(0)},
$$
then we have
$$
\lVert A_+-A\rVert\leqslant \lVert A\rVert \lVert Id-e^{\hat{f}^{re}(0)} \rVert \leqslant 2\lVert A\rVert\epsilon.
$$

  \textbf{Resonant case:}
 In fact, we only need to consider the case in which $A$ is elliptic with eigenvalues $\{e^{i\rho},e^{-i\rho}\}$ for $\rho\in \R\backslash\{0\}$ since if $\rho\in i\R$, then the non-resonant condition is always satisfied due to the Diophantine condition on $\alpha$ and then it actually belongs to the non-resonant case.

\begin{Claim} $n_\ast$ is the unique resonant site with
$$
0<\lvert n_\ast\rvert \leqslant N=\frac{2}{r-r'} \lvert \ln \epsilon \rvert.
$$
\end{Claim}
\begin{pf}Indeed, if there exists  $n_{\ast}^{'}\neq n_\ast$ satisfying $|2\rho- \la n_{\ast}^{'},\alpha\ra|< \epsilon^{\sigma}$, then by the Diophantine condition of $\alpha$, we have
$$
\frac{\kappa}{\lvert n_{\ast}^{'}-n_\ast\rvert^{\tau}}\leqslant \lvert \la n_{\ast}^{'}-n_\ast,\alpha\ra\rvert< 2\epsilon^{\sigma},
$$
which implies that
$\lvert n_{\ast}^{'} \rvert>2^{-\frac{1}{\tau}}\kappa^{\frac{1}{\tau}}\epsilon^{-\frac{\sigma}{\tau}}-N> 2N^2.$\end{pf}

Since we have \begin{equation}\label{reso}
\lvert 2\rho- \la n_\ast,\alpha\ra \rvert< \epsilon^{\sigma},
\end{equation}
 the smallness condition on $\epsilon$  implies that

$$\lvert \ln \epsilon\lvert ^\tau \epsilon^\sigma \leqslant \frac{\kappa (r-r')^\tau}{2^{\tau+1}}.$$

\noindent Thus

$$\frac{\kappa}{\lvert n_\ast\lvert^\tau} \leqslant \lvert \la n_\ast,\alpha\ra \lvert\leqslant \epsilon^\sigma +2\lvert \rho\lvert \leqslant \frac{\kappa}{2\lvert n_\ast\lvert ^\tau}+2\lvert \rho\lvert ,$$

\noindent which implies that

$$
\lvert \rho \rvert\geqslant \frac{\kappa}{4\lvert n_\ast\rvert^{\tau}}.
$$
Then by Lemma 8.1 of Hou-You \cite{houyou}, one can find $P\in SU(1,1)$ with
$$
\lVert P \rVert \leqslant 2(\frac{\lVert A \rVert}{\lvert \rho \rvert})^{\frac{1}{2}}\leqslant 4(\frac{\lVert A\rVert}{\kappa})^{\frac{1}{2}}\lvert n_\ast \rvert^{\frac{\tau}{2}},
$$
such that
$$
PAP^{-1}=\begin{pmatrix} e^{i\rho} & 0\\ 0 & e^{-i\rho} \end{pmatrix}=A'.
$$
Denote $g=PfP^{-1}$, by $(\ref{estf})$ we have:
\begin{eqnarray}
\label{esti-p}\lVert P \rVert &\leqslant& 4(\frac{\lVert A\rVert}{\kappa})^{\frac{1}{2}}\lvert N \rvert^{\frac{\tau}{2}}\leqslant 4(\frac{\lVert A\rVert}{\kappa})^{\frac{1}{2}}(\frac{2}{r-r'} \lvert \ln \epsilon \rvert)^{\frac{\tau}{2}},\\
\label{esti-g}\lvert g \rvert_r &\leqslant& \lVert P \rVert^2\lvert f\rvert_r \leqslant  \frac{2^{4+\tau}\lVert A\rVert\lvert \ln\epsilon\rvert^{\tau}}{\kappa(r-r')^{\tau}}\times\epsilon:=\epsilon'.
\end{eqnarray}

Now we define
\begin{flalign*}
&\Lambda_1(\epsilon^{\sigma})=\{n\in\Z^{d}: \lvert \la n,\alpha\ra\rvert \geqslant \epsilon^{\sigma}\},\\
&\Lambda_2(\epsilon^{\sigma})=\{n\in\Z^{d}: \lvert 2\rho-\la n,\alpha\ra\rvert \geqslant \epsilon^{\sigma}\}.
\end{flalign*}
For $\eta=\epsilon^{\sigma}$, we define the decomposition $\mathcal{B}_r=\mathcal{B}_r^{nre}(\epsilon^{\sigma}) \bigoplus\mathcal{B}_r^{re}(\epsilon^{\sigma})$ as in $(\ref{space})$ with $A$ substituted by $A'$. Direct computation shows that any $Y\in \mathcal{B}_r^{nre}(\epsilon^{\sigma})$ takes the precise form:
\begin{small}
\begin{equation}
\begin{split}
Y(\theta)=&\sum_{n\in \Lambda_1(\epsilon^{\sigma})}\begin{pmatrix} i\hat{t}(n) & 0\\ 0 & -i\hat{t}(n) \end{pmatrix} e^{i\la n,\theta\ra}+\\
&\sum_{n\in \Lambda_2(\epsilon^{\sigma})}\begin{pmatrix} 0 & \hat{v}(n)e^{i\la n,\theta\ra}\\ \overline{\hat{v}(n)}e^{-i\la n,\theta\ra} & 0 \end{pmatrix}.
\end{split}
\end{equation}
\end{small}

Since $\epsilon^{\sigma}\geqslant 13\lVert A'\rVert^2 (\epsilon') ^{\frac{1}{2}}$, we can apply Lemma \ref{lem2} to remove all the non-resonant terms of $g$, which means there exist $Y\in \mathcal{B}_r$ and $g^{re}\in \mathcal{B}_r^{re}(\eta)$ such that
$$
e^{Y(\theta+\alpha)}(A'e^{g(\theta)})e^{-Y(\theta)}=A'e^{g^{re}(\theta)},
$$
with $\lvert Y \rvert_r\leqslant (\epsilon')^{\frac{1}{2}}$ and $\lvert g^{re}\rvert_r\leqslant 2\epsilon'$.

Combining with the Diophantine condition on the frequency $\alpha$ and the Claim, we have:
\begin{flalign*}
&\{\Z^{d}\backslash\Lambda_1(\epsilon^{\sigma})\}\cap\{n\in \Z^{d}:\lvert n \rvert\leqslant \kappa^{\frac{1}{\tau}}\epsilon^{-\frac{\sigma}{\tau}}\}=\{0\},\\
&\{\Z^{d}\backslash\Lambda_2(\epsilon^{\sigma})\}\cap\{n\in \Z^{d}:\lvert n \rvert\leqslant 2^{-\frac{1}{\tau}}\kappa^{\frac{1}{\tau}}\epsilon^{-\frac{\sigma}{\tau}}-N\}=\{n_\ast\}.
\end{flalign*}
Let $N':=2^{-\frac{1}{\tau}}\kappa^{\frac{1}{\tau}}\epsilon^{-\frac{\sigma}{\tau}}-N$, then we can rewrite $g^{re}(\theta)$ as
\begin{flalign*}
g^{re}(\theta)&=g^{re}_0+g^{re}_1(\theta)+g^{re}_2(\theta)\\
&=\begin{pmatrix} i\hat{t}(0) & 0 \\ 0 & -i\hat{t}(0) \end{pmatrix}+\begin{pmatrix} 0 & \hat{v}(n_\ast)e^{i\la n_\ast,\theta\ra} \\ \overline{\hat{v}(n_\ast)}e^{-i\la n_\ast,\theta\ra} & 0 \end{pmatrix}\\
&+\sum_{\lvert n \rvert>N'}\hat{g}^{re}(n)e^{i\la n,\theta\ra}.
\end{flalign*}

Define the $4\pi\Z^d$-periodic rotation $Q(\theta)$ as below:
$$
Q(\theta)=\begin{pmatrix} e^{-\frac{\la n_\ast,\theta\ra}{2}i} & 0\\ 0 & e^{\frac{\la n_\ast,\theta\ra}{2}i} \end{pmatrix}.
$$
So we have \begin{equation}\label{esti-Q}\lvert Q(\theta)\rvert_{r'}\leqslant e^{\frac{1}{2}Nr'}\leqslant \epsilon^{\frac{-r'}{r-r'}}.\end{equation} One can also show that
$$
Q(\theta+\alpha)(A'e^{g^{re}(\theta)})Q^{-1}(\theta)=\tilde{A}e^{\tilde{g}(\theta)},
$$
where
\begin{equation}\label{tildea}
\tilde{A}=Q(\theta+\alpha)A'Q^{-1}(\theta)=\begin{pmatrix} e^{i(\rho-\frac{\la n_\ast,\alpha\ra}{2})} & 0\\ 0 & e^{-i(\rho-\frac{\la n_\ast,\alpha\ra}{2})} \end{pmatrix}
\end{equation}
and
$$\tilde{g}(\theta)=Qg^{re}(\theta)Q^{-1}=Qg^{re}_0Q^{-1}+Qg^{re}_1(\theta)Q^{-1}+Qg^{re}_2(\theta)Q^{-1}.$$
Moreover,
\begin{flalign}
\label{ell1}Qg^{re}_0Q^{-1} &=g^{re}_0 = \begin{pmatrix} i\hat{t}(0) & 0 \\ 0 & -i\hat{t}(0) \end{pmatrix}
\in su(1,1), \\
\label{ell2}Qg^{re}_1(\theta)Q^{-1}&=\begin{pmatrix} 0 & \hat{v}(n_\ast) \\ \overline{\hat{v}(n_\ast)} & 0 \end{pmatrix} \in su(1,1).
\end{flalign}

Now we return back from $su(1,1)$ to $sl(2,\R)$. Denote
\begin{flalign}
\label{ell}L &=M^{-1}(Qg^{re}_0Q^{-1}+Qg^{re}_1(\theta)Q^{-1})M,\\
F &=M^{-1}Qg^{re}_2(\theta)Q^{-1}M,\\
B &=M^{-1}(Q\circ e^Y \circ P) M,\\
\label{tildeaa}\tilde{A}^{'} &=M^{-1}\tilde{A}M,
\end{flalign}
then we have:
\begin{equation}\label{con1}
B(\theta+\alpha)(Ae^{f(\theta)})B^{-1}(\theta)=\tilde{A}^{'}e^{L+F(\theta)}.
\end{equation}
By $(\ref{esti-p})$ and $(\ref{esti-Q})$, we have the following estimates:
\small\begin{eqnarray}
 \lVert B\rVert_0 &\leqslant&  |e^{Y}|_r \lVert P\rVert \leqslant 8(\frac{\lVert A\rVert}{\kappa})^{\frac{1}{2}}(\frac{2}{r-r'} \lvert \ln \epsilon \rvert)^{\frac{\tau}{2}},\\
\lvert B\rvert_{r'} &\leqslant& 8(\frac{\lVert A\rVert}{\kappa})^{\frac{1}{2}}(\frac{2}{r-r'} \lvert \ln \epsilon \rvert)^{\frac{\tau}{2}}\times\epsilon^{\frac{-r'}{r-r'}},\\
\label{D}\lVert L \rVert &\leqslant& \lVert Qg^{re}_0Q^{-1}\rVert + \lVert Qg^{re}_1(\theta)Q^{-1}\rVert \leqslant \epsilon'+\epsilon' e^{-\lvert n_{\ast}\rvert r},\\
\label{F}|F|_{r'} &\leqslant& \lvert Qg^{re}_2(\theta)Q^{-1}\rvert_{r'} \leqslant \frac{2^{4+\tau}\lVert A\rVert\lvert \ln\epsilon\rvert^{\tau}}{\kappa(r-r')^{\tau}}\epsilon e^{-N'(r-r')}(N')^de^{Nr'}.
\end{eqnarray}

By $(\ref{D})$ and $(\ref{F})$, direct computation shows that
\begin{equation}\label{impl}
e^{L+F(\theta)}=e^L+\mathcal{O}(F(\theta))=e^L(Id+e^{-L}\mathcal{O}(F(\theta)))=e^L e^{f_+{(\theta)}}.
\end{equation}
It immediately implies that
$$
\lvert f_+{(\theta)}\rvert_{r'}\leqslant 2 |F(\theta)|_{r'}\leqslant \frac{2^{5+\tau}\lVert A\rVert\lvert \ln\epsilon\rvert^{\tau}}{\kappa(r-r')^{\tau}}\epsilon e^{-N'(r-r')}(N')^de^{Nr'}\ll \epsilon^{100}.
$$
Thus we can rewrite  $(\ref{con1})$ as
$$
B(\theta+\alpha)(Ae^{f(\theta)})B^{-1}(\theta)=A_{+}e^{f_+(\theta)},
$$
with
\begin{equation}\label{constm}
A_+=\tilde{A}^{'}e^L=e^{A''}, \ \ A''\in sl(2,\R).
\end{equation}

Now recall that Baker-Campbell-Hausdorff Formula says that
\begin{equation}
\ln(e^X e^Y)=X+Y+\frac{1}{2}[X,Y]+\frac{1}{12}([X,[X,Y]+[Y,[Y,X]])+\cdots,
\end{equation}
where $[X,Y]=XY-YX$ denotes the Lie Bracket and $\cdots$ denotes the sum of higher order terms. Using this formula and by a simple calculation, $(\ref{constm})$ gives
$$
MA''M^{-1}=\begin{pmatrix} it & v\\ \bar{v} & -it \end{pmatrix}
$$
where
$$
t=\rho-\frac{\la n_\ast,\alpha\ra}{2}+\hat{t}(0)+ {\rm higher\, order\, terms}
$$
and
$$
v=\hat{v}(n_\ast)+ {\rm higher\, order\, terms}.
$$
By $(\ref{reso})$ and $(\ref{D})$, we obtain $\lvert t\rvert \leqslant \epsilon^{\sigma}$ and
$$
\lvert v \rvert\leqslant \frac{2^{4+\tau}\lVert A\rVert\lvert \ln\epsilon\rvert^{\tau}}{\kappa(r-r')^{\tau}}\epsilon e^{-\lvert n_{\ast}\rvert r}.
$$
Finally, the following estimate is straightforward:
\begin{equation}\label{constant}
\lVert A'' \rVert\leqslant 2(\lvert \rho-\frac{\la n_\ast,\alpha\ra}{2}\rvert+\lVert Qg^{re}_0Q^{-1}\rVert+\lVert Qg^{re}_1(\theta)Q^{-1}\rVert)\leqslant 2\epsilon^{\sigma}.
\end{equation}
This finishes the proof of Proposition $\ref{prop1}$.
\end{pf}

\begin{Remark}
The special structure of $A_+$ can give us a precise estimate of the upper triangular element of the parabolic constant matrix when the rotation number of the initial system is rational with respect to $\alpha$, see \cite{LYZZ}. However, it does not work for general case. Instead, $A_+=e^{A''}$ with $\lVert A''\rVert \leqslant 2\epsilon^{\sigma}$ will help.
\end{Remark}

\begin{Remark}\label{rem3.3}
The limitation $\sigma<\frac{1}{6}$ not only makes Lemma \ref{lem2} applicable but also ensures that our conjugation map $B$ takes value in $SL(2,\R)$. The choice of $D$ being larger than $\frac{2}{\sigma}$ is necessary for us to guarantee the arbitrariness of $r'\in(0,r)$ and the separation of resonant steps $($see Claim \ref{cl2}$)$.
\end{Remark}

\subsection{$C^k$ Almost Reducibility}

Let $(f_{j})_{j\geqslant 1}$, $f_{j}\in C_{\frac{1}{j}}^{\omega}(\T^{d},sl(2,\R))$ be the analytic sequence approximating $f\in C^k(\T^{d},sl(2,\R))$ which satisfies $(\ref{2.1})$.

For $0<r'<r$, denote
\begin{equation}
\epsilon_0^{'}(r,r')=\frac{c}{(2\lVert A\rVert)^{\tilde{D}}}(r-r')^{D\tau},
\end{equation}
where $c,D,\tilde{D},\tau$ are defined in Theorem \ref{prop1}.

For $m\in\Z^+$, we define
\begin{equation}
\epsilon_m=\frac{c}{(2\lVert A\rVert)^{\tilde{D}}m^{D\tau+\frac{1}{2}}}.
\end{equation}

Then for any $0<s\leqslant\frac{1}{6D\tau+3}$ fixed, there exists $m_0$ such that for any $m\geqslant m_0$ we have both
\begin{equation}\label{pickm0}
\epsilon_m\leqslant \epsilon_0^{'}(\frac{1}{m},\frac{1}{m^{1+s}}),
\end{equation}
and
$$
\frac{1}{m^s-1}\leqslant \frac{s}{4}.
$$
We will start from $M> \max\{\frac{(2\lVert A\rVert)^{\tilde{D}}}{c},m_0\}$, $M\in\N^+$. Denote $l_j=M^{(1+s)^{j-1}}$, $j\in \N^+$. In case that $l_j$ is not an integer, we just pick $[l_j]+1$ instead of $l_j$.

Now, denote by $\Omega=\{l_{n_1},l_{n_2},l_{n_3},\cdots\}$ the sequence of all resonant steps. That is, $l_{n_j}$-th step is obtained by resonant case. Using analytic approximation $(\ref{2.1})$ and Theorem \ref{prop1}, we have the following Proposition concerning each $(\alpha, Ae^{f_{l_j}(\theta)})$.

\begin{Proposition}\label{pro33}
Let $\alpha\in {\rm DC}(\kappa,\tau)$, $\sigma<\frac{1}{6}$. Assume that $A\in SL(2,\R)$, $f\in C^k(\T^d, sl(2,\R))$ with $k>(D+2)\tau+2$ and $\{f_j\}_{j\geqslant 1}$ are defined above. There exists $\epsilon_0=\epsilon_0(\kappa,\tau,d,k, \lVert A\rVert,\sigma)$ such that if $\lVert f\rVert_k \leqslant \epsilon_0$, then there exist $B_{l_j}\in C^{\omega}_{\frac{1}{l_{j+1}}}(2\T^d, SL(2,\R))$, $A_{l_j}\in SL(2,\R)$ and $f^{'}_{l_j}\in C^{\omega}_{\frac{1}{l_{j+1}}}(2\T^d, sl(2,\R))$ such that
$$
B_{l_j}(\theta+\alpha)(Ae^{f_{l_j}(\theta)})B^{-1}_{l_j}(\theta)=A_{l_j}e^{f_{l_j}^{'}(\theta)},
$$
with estimates

\begin{flalign}\label{estimBlj}\lvert B_{l_j}(\theta)\rvert_{\frac{1}{l_{j+1}}}&\leqslant 64(\frac{\lVert A\rVert}{\kappa})(\frac{2}{\frac{1}{l_j}-\frac{1}{l_{j+1}}} \lvert \ln \epsilon_{l_j} \rvert)^{\tau}\times{\epsilon_{l_j}}^{-\frac{\frac{2}{l_{j+1}}}{\frac{1}{l_j}-\frac{1}{l_{j+1}}}} \leqslant \epsilon_{l_{j}}^{-\frac{\sigma}{2}-s},\\
\label{estimBlj2}\lVert B_{l_j}(\theta)\rVert_0&\leqslant 64(\frac{\lVert A\rVert}{\kappa})(\frac{2}{\frac{1}{l_j}-\frac{1}{l_{j+1}}} \lvert \ln \epsilon_{l_j} \rvert)^{\tau} \leqslant \epsilon_{l_{j}}^{-\frac{\sigma}{2}},\\
\label{estimflj}\lvert f_{l_j}^{'}(\theta)\rvert_{\frac{1}{l_{j+1}}}&\leqslant \epsilon_{l_{j}}^{3-\sigma},\ \ \lVert A_{l_j}\rVert\leqslant 2\lVert A\rVert.
\end{flalign}

Moreover, there exists unitary matrices $U_j\in SL(2,\C)$ such that
$$
U_jA_{l_j}U_j^{-1}=\begin{pmatrix} e^{\gamma_j} & c_j\\ 0 & e^{-\gamma_j} \end{pmatrix}
$$
and
\begin{equation}
\label{estisharp}\lVert B_{l_j}(\theta)\rVert_0^2\lvert c_j\rvert \leqslant 8\lVert A\rVert,
\end{equation}
with $\gamma_j\in i\R\cup\R$ and $c_j\in \C$.

\end{Proposition}

\begin{pf} \textbf{First step:} Assume that
$$
C'\lVert f(\theta) \rVert_k\leqslant \frac{c}{(2\lVert A\rVert)^{\tilde{D}}l_1^{D\tau+\frac{1}{2}}},
$$
then by $(\ref{2.1})$ and $(\ref{pickm0})$ we have
$$
\lvert f_{l_{1}}(\theta)\rvert_{\frac{1}{l_{1}}}\leqslant \epsilon_{l_{1}}\leqslant \epsilon_0^{'}(\frac{1}{l_{1}},\frac{1}{l_{2}}).
$$
Apply Theorem \ref{prop1}, we can find $B_{l_1}\in C^{\omega}_{\frac{1}{l_2}}(2\T^{d},SL(2,\R))$, $A_{l_1}\in SL(2,\R)$ and $f_{l_1}^{'}\in C^{\omega}_{\frac{1}{l_2}}(\T^{d},sl(2,\R))$ such that
$$
B_{l_1}(\theta+\alpha)(Ae^{f_{l_1}(\theta)})B^{-1}_{l_1}(\theta)=A_{l_1}e^{f_{l_1}^{'}(\theta)}.
$$
More precisely, we have two different cases:
\begin{itemize}
\item (Non-resonant case)
$$
\lvert B_{l_1}\rvert_{\frac{1}{l_2}}\leqslant 1+\epsilon_{l_{1}}^{\frac{1}{2}}, \ \ \lvert f_{l_1}^{'}\rvert_{\frac{1}{l_2}}\leqslant \epsilon_{l_{1}}^{3-\sigma},
$$
and
$$
\lVert A_{l_1}-A\rVert \leqslant 2\lVert A\rVert\epsilon_{l_{1}}.
$$
\item (Resonant case)
\begin{flalign*}
\lvert B_{l_1}\rvert_{\frac{1}{l_2}} &\leqslant 8(\frac{\lVert A\rVert}{\kappa})^{\frac{1}{2}}(\frac{2}{\frac{1}{l_1}-\frac{1}{l_2}} \lvert \ln \epsilon_{l_1} \rvert)^{\frac{\tau}{2}}\times{\epsilon_{l_1}}^{-\frac{\frac{1}{l_2}}{\frac{1}{l_1}-\frac{1}{l_2}}},\\
\lvert B_{l_1}\rvert_0 &\leqslant 8(\frac{\lVert A\rVert}{\kappa})^{\frac{1}{2}}(\frac{2}{\frac{1}{l_1}-\frac{1}{l_2}} \lvert \ln \epsilon_{l_1} \rvert)^{\frac{\tau}{2}}, \ \ \lvert f_{l_1}^{'}\rvert_{\frac{1}{l_2}}\leqslant \epsilon_{l_{1}}^{100}.
\end{flalign*}
Moreover, $A_{l_1}=e^{A_{l_1}''}$ with $\lVert A_{l_1}''\rVert \leqslant 2\epsilon_{l_{1}}^{\sigma}$.
\end{itemize}
In both cases, it is clear that $(\ref{estimBlj})$, $(\ref{estimBlj2})$, $(\ref{estimflj})$ and $(\ref{estisharp})$ are fulfilled.

Note that the first step is a little special as it does not involve the composition of conjugation maps. Therefore, in order to provide a more explicit iteration process, let us show one more step before induction.

\textbf{Second step:}
We have
$$
B_{l_1}(\theta+\alpha)(Ae^{f_{l_{2}}})B^{-1}_{l_1}(\theta)=A_{l_1}e^{f_{l_1}^{'}}+B_{l_1}(\theta+\alpha)(Ae^{f_{l_{2}}}-Ae^{f_{l_{1}}})B^{-1}_{l_1}(\theta).
$$
We can rewrite that $$
A_{l_1}e^{f_{l_1}^{'}(\theta)}+B_{l_1}(\theta+\alpha)(Ae^{f_{l_{2}}(\theta)}-Ae^{f_{l_{1}}(\theta)})B^{-1}_{l_1}(\theta)=A_{l_1}e^{\widetilde{f_{l_1}}(\theta)}.
$$

We pick $k>(D+2)\tau+2\geqslant (1+3s)(D\tau+\frac{1}{2})+2\tau+1.$

If the previous step is non-resonant,
\begin{flalign*}
\lvert \widetilde{f_{l_1}}(\theta)\rvert_{\frac{1}{l_{2}}}\leqslant & \lvert f_{l_1}^{'}(\theta)\rvert_{\frac{1}{l_{2}}}+\lVert A_{l_1}^{-1}\rVert\lvert B_{l_1}(\theta+\alpha)(Ae^{f_{l_{2}}(\theta)}-Ae^{f_{l_{1}}(\theta)})B^{-1}_{l_1}(\theta)\rvert_{\frac{1}{l_{2}}}\\
\leqslant & 4\epsilon_{l_{1}}^{3-\sigma}+ 2\lVert A\rVert^2\times2\times\frac{c}{(2\lVert A\rVert)^{\tilde{D}}l_1^{D\tau+\frac{1}{2}}l_{1}^{k-1}}\\
\leqslant & \epsilon_{l_2}\\
\leqslant & \epsilon_0^{'}(\frac{1}{l_{2}},\frac{1}{l_{3}}).
\end{flalign*}

If the previous step is resonant,
\begin{flalign*}
\lvert \widetilde{f_{l_1}}(\theta)\rvert_{\frac{1}{l_{2}}}\leqslant & \lvert f_{l_1}^{'}(\theta)\rvert_{\frac{1}{l_{2}}}+\lVert A_{l_1}^{-1}\rVert\lvert B_{l_1}(\theta+\alpha)(Ae^{f_{l_{2}}(\theta)}-Ae^{f_{l_{1}}(\theta)})B^{-1}_{l_1}(\theta)\rvert_{\frac{1}{l_{2}}}\\
\leqslant & \epsilon_{l_{1}}^{100}+ 128(\frac{\lVert A\rVert^3}{\kappa})(\frac{2}{\frac{1}{l_1}-\frac{1}{l_2}}\lvert \ln \epsilon_{l_1} \rvert)^{\tau}{\epsilon_{l_1}}^{- \frac{\frac{2}{l_2}}{\frac{1}{l_1}-\frac{1}{l_2}}}\frac{c}{(2\lVert A\rVert)^{\tilde{D}}l_1^{D\tau+\frac{1}{2}}l_{1}^{k-1}}\\
\leqslant & \epsilon_{l_2}\\
\leqslant & \epsilon_0^{'}(\frac{1}{l_{2}},\frac{1}{l_{3}}).
\end{flalign*}

Now for $(\alpha,A_{l_1}e^{\widetilde{f_{l_1}}(\theta)})$, we can apply Proposition \ref{prop1} again to get $\tilde{B}_{l_1}\in C^{\omega}_{\frac{1}{l_{3}}}(2\T^{d},$
$SL(2,\R))$, $A_{l_{2}}\in SL(2,\R)$ and $f_{l_{2}}^{'}\in C^{\omega}_{\frac{1}{l_{3}}}(\T^{d},sl(2,\R))$ such that
\begin{equation}\label{second}
\tilde{B}_{l_1}(\theta+\alpha)(A_{l_1}e^{\widetilde{f_{l_1}}(\theta)})\tilde{B}_{l_1}^{-1}(\theta)=A_{l_{2}}e^{f_{l_{2}}^{'}(\theta)},
\end{equation}
which gives
\begin{equation}
B_{l_2}(\theta+\alpha)(Ae^{f_{l_2}(\theta)})B^{-1}_{l_2}(\theta)=A_{l_2}e^{f_{l_2}^{'}(\theta)}.
\end{equation}
Here $B_{l_2}(\theta)=\tilde{B}_{l_1}(\theta)\circ B_{l_1}(\theta)$.

Before giving the precise estimates of each term, let us first introduce a Claim showing that all the resonant steps are separated.
\begin{Claim}\label{cl2}
We have $\epsilon_{l_{n_{j+1}}}<\epsilon^2_{l_{n_{j}}}, \forall j\in \Z^+$.
\end{Claim}
\begin{pf}
It is enough to prove for $j=1$. By definition, we have $l_{n_{1}}$-th step is obtained by resonant case. Thus Theorem \ref{prop1} implies $\rho(A_{l_{n_{1}}})\leqslant 2\epsilon_{l_{n_{1}}}^{\sigma}$ and $\rho(A_{l_{n_{2}-1}})\leqslant 4\epsilon_{l_{n_{1}}}^{\sigma}$ since every step between $l_{n_{1}}$ and $l_{n_{2}}$ is non-resonant.

By the resonant condition of $l_{n_{2}}$-th step, there exists a unique $n$ with $0<\lvert n\rvert \leqslant N_{l_{n_{2}}}$ such that
\begin{equation}\label{resonant}
\lvert 2\rho(A_{l_{n_{2}-1}})- \la n,\alpha \ra\rvert\leqslant \epsilon_{l_{n_{2}}}^{\sigma}.
\end{equation}
However, by the Diophantine condition of $\alpha$ and condition $(\ref{estf})$, we have
$$
\lvert \la n,\alpha\ra\rvert \geqslant \frac{\kappa}{\lvert N_{l_{n_{2}}}\rvert^{\tau}}\geqslant 10\epsilon_{l_{n_{2}}}^{\frac{\sigma}{2}}.
$$
If $\epsilon_{l_{n_{2}}}\geqslant\epsilon^2_{l_{n_{1}}}$, then $(\ref{resonant})$ yields contradiction.
\end{pf}

Now by Claim \ref{cl2}, we have (in the worst situation)
\begin{flalign*}
\lvert B_{l_2}(\theta)\rvert_{\frac{1}{l_{3}}} & \leqslant 64(\frac{\lVert A\rVert}{\kappa})(\frac{2}{\frac{1}{l_2}-\frac{1}{l_{3}}} \lvert \ln \epsilon_{l_2} \rvert)^{\tau}\times{\epsilon_{l_2}}^{-\frac{\frac{2}{l_{3}}}{\frac{1}{l_2}-\frac{1}{l_{3}}}}\leqslant\epsilon_{l_2}^{-\frac{\sigma}{2}-s}\\
\lVert B_{l_2}(\theta)\rVert_0 & \leqslant 64(\frac{\lVert A\rVert}{\kappa})(\frac{2}{\frac{1}{l_2}-\frac{1}{l_{3}}} \lvert \ln \epsilon_{l_2} \rvert)^{\tau}\leqslant \epsilon_{l_2}^{-\frac{\sigma}{2}}\\
\lvert f_{l_{2}}^{'}(\theta)\rvert_{\frac{1}{l_{3}}} & \leqslant \epsilon_{l_2}^{3-\sigma}, \ \ \lVert A_{l_j}\rVert \leqslant 2\lVert A\rVert.
\end{flalign*}
Thus $(\ref{estimBlj})$, $(\ref{estimBlj2})$ and $(\ref{estimflj})$ are fulfilled again. For $(\ref{estisharp})$, all three cases (1. there is no resonance, 2. the first step is resonant, 3. the second step is resonant) satisfy it. The discussion is similar to that in the induction step below, so we omit it here for simplicity.

\textbf{Induction step:} Assume that for $l_n,n\leqslant \tilde{n}$, we already have $(\ref{estisharp})$ and
\begin{equation}
\label{estiln}B_{l_n}(\theta+\alpha)(Ae^{f_{l_n}(\theta)})B^{-1}_{l_n}(\theta)=A_{l_n}e^{f_{l_n}^{'}(\theta)},
\end{equation}
with
\begin{equation}
\label{estima}\lvert B_{l_n}(\theta)\rvert_{\frac{1}{l_{n+1}}} \leqslant 64(\frac{\lVert A\rVert}{\kappa})(\frac{2}{\frac{1}{l_n}-\frac{1}{l_{n+1}}} \lvert \ln \epsilon_{l_n} \rvert)^{\tau}{\epsilon_{l_n}}^{-\frac{\frac{2}{l_{n+1}}}{\frac{1}{l_n}-\frac{1}{l_{n+1}}}} \triangleq \zeta_n^2 \leqslant \epsilon_{l_n}^{-\frac{\sigma}{2}-s},
\end{equation}
\begin{equation}\label{estimb}
\lVert B_{l_n}(\theta)\rVert_0 \leqslant 64(\frac{\lVert A\rVert}{\kappa})(\frac{2}{\frac{1}{l_n}-\frac{1}{l_{n+1}}} \lvert \ln \epsilon_{l_n} \rvert)^{\tau} \leqslant \epsilon_{l_n}^{-\frac{\sigma}{2}},\ \ \lvert f_{l_n}^{'}(\theta)\rvert_{\frac{1}{l_{n+1}}} \leqslant \epsilon_{l_n}^{3-\sigma},
\end{equation}
and
\begin{equation}
\label{esta}\lVert A_{l_n}\rVert \leqslant 2\lVert A\rVert.
\end{equation}
Moreover, if the $n$-th step is obtained by the resonant case, we have
\begin{equation}
\label{estres}A_{l_n}=e^{A_{l_n}''}, \ \ \lVert A_{l_n}''\rVert <2\epsilon_{l_{n}}^{\sigma}.
\end{equation}
If the $n$-th step is obtained by the non-resonant case, we have
\begin{equation}
\label{estnon}\lVert A_{l_n}-A_{l_{n-1}}\rVert\leqslant 2\lVert A_{l_{n-1}}\rVert \epsilon_{l_{n}}.
\end{equation}
and
\begin{equation}
\label{estb}\lVert B_{l_n}\rVert_{\frac{1}{l_{n+1}}}\leqslant (1+\epsilon_{l_{n}}^{\frac{1}{2}})\lVert B_{l_{n-1}}\rVert_{\frac{1}{l_{n}}},\ \ \lVert B_{l_n}\rVert_{0}\leqslant (1+\epsilon_{l_{n}}^{\frac{1}{2}})\lVert B_{l_{n-1}}\rVert_{0}.
\end{equation}
Now by $(\ref{estiln})$, for $l_{n+1}, n=\tilde{n}$, we have
$$
B_{l_n}(\theta+\alpha)(Ae^{f_{l_{n+1}}})B^{-1}_{l_n}(\theta)=A_{l_n}e^{f_{l_n}^{'}}+B_{l_n}(\theta+\alpha)(Ae^{f_{l_{n+1}}}-Ae^{f_{l_{n}}})B^{-1}_{l_n}(\theta).
$$
In $(\ref{2.1})$, a simple summation implies
\begin{equation}\label{ov}
\lvert f_{l_{n+1}}(\theta)-f_{l_{n}}(\theta)\rvert_{\frac{1}{l_{n+1}}}\leqslant \frac{c}{(2\lVert A\rVert)^{\tilde{D}} l_1^{D\tau+\frac{1}{2}}l_{n}^{k-1}},
\end{equation}
Moreover, $(\ref{2.1})$ also gives us
\begin{equation}\label{ok}
\lvert f_{l_{n+1}}(\theta)\rvert_{\frac{1}{l_{n+1}}} +\lvert f_{l_{n}}(\theta)\rvert_{\frac{1}{l_{n+1}}}\leqslant \frac{2c}{(2\lVert A\rVert)^{\tilde{D}} M^{D\tau+\frac{1}{2}}}.
\end{equation}
Thus if we rewrite that $$A_{l_n}e^{f_{l_n}^{'}(\theta)}+B_{l_n}(\theta+\alpha)(Ae^{f_{l_{n+1}}(\theta)}-Ae^{f_{l_{n}}(\theta)})B^{-1}_{l_n}(\theta)=A_{l_n}e^{\widetilde{f_{l_n}}(\theta)},$$
by $(\ref{estima})$, $(\ref{estimb})$, $(\ref{esta})$, $(\ref{ov})$ and $(\ref{ok})$ we obtain
\begin{small}
\begin{flalign*}
\lvert \widetilde{f_{l_n}}(\theta)\rvert_{\frac{1}{l_{n+1}}}\leqslant & \lvert f_{l_n}^{'}(\theta)\rvert_{\frac{1}{l_{n+1}}}+\lVert A_{l_n}^{-1}\rVert\lvert B_{l_n}(\theta+\alpha)(Ae^{f_{l_{n+1}}(\theta)}-Ae^{f_{l_{n}}(\theta)})B^{-1}_{l_n}(\theta)\rvert_{\frac{1}{l_{n+1}}}\\
\leqslant & \epsilon_{l_{n}}^{3-\sigma}+ (\frac{\lVert A\rVert^4}{\kappa^2})(\frac{2}{\frac{1}{l_n}-\frac{1}{l_{n+1}}}\lvert \ln \epsilon_{l_n} \rvert)^{2\tau}{\epsilon_{l_n}}^{- \frac{\frac{4}{l_{n+1}}}{\frac{1}{l_n}-\frac{1}{l_{n+1}}}}\frac{2\times 64\times 64 c}{(2\lVert A\rVert)^{\tilde{D}}l_1^{D\tau+\frac{1}{2}}l_{n}^{k-1}}\\
\leqslant & \epsilon_{l_{n+1}}\\
\leqslant & \epsilon_0^{'}(\frac{1}{l_{n+1}},\frac{1}{l_{n+2}}).
\end{flalign*}
\end{small}
Now for $(\alpha,A_{l_n}e^{\widetilde{f_{l_n}}(\theta)})$, we can apply Proposition \ref{prop1} again to get $\tilde{B}_{l_n}\in C^{\omega}_{\frac{1}{l_{n+2}}}(2\T^{d},$
$SL(2,\R))$, $A_{l_{n+1}}\in SL(2,\R)$ and $f_{l_{n+1}}^{'}\in C^{\omega}_{\frac{1}{l_{n+2}}}(\T^{d},sl(2,\R))$ such that
$$
\tilde{B}_{l_n}(\theta+\alpha)(A_{l_n}e^{\widetilde{f_{l_n}}(\theta)})\tilde{B}_{l_n}^{-1}(\theta)=A_{l_{n+1}}e^{f_{l_{n+1}}^{'}(\theta)},
$$
with
$$\lvert f_{l_{n+1}}^{'}(\theta)\rvert_{\frac{1}{l_{n+2}}}\leqslant \epsilon_{l_{n+1}}^{3-\sigma}.$$

Denote $B_{l_{n+1}}:=\tilde{B}_{l_n}B_{l_{n}} \in C^{\omega}_{\frac{1}{l_{n+2}}}(2\T^{d},SL(2,\R))$, by Claim \ref{cl2}, we have (in the worst situation)
\begin{flalign*}
\lvert B_{l_{n+1}}(\theta)\rvert_{\frac{1}{l_{n+2}}}&\leqslant \zeta_{n+1}^{1+\frac{1}{2}+(\frac{1}{2})^2+(\frac{1}{2})^3+\cdots+(\frac{1}{2})^n+\cdots}\\
&\leqslant \zeta_{n+1}^2\\
&\leqslant 64(\frac{\lVert A\rVert}{\kappa})(\frac{2}{\frac{1}{l_{n+1}}-\frac{1}{l_{n+2}}} \lvert \ln \epsilon_{l_{n+1}} \rvert)^{\tau}\times{\epsilon_{l_{n+1}}}^{-\frac{\frac{2}{l_{n+2}}}{\frac{1}{l_{n+1}}-\frac{1}{l_{n+2}}}} \\
&\leqslant \epsilon_{l_{n+1}}^{-\frac{\sigma}{2}-s}.
\end{flalign*}
and
$$
\lVert B_{l_{n+1}}(\theta) \rVert_0\leqslant 64(\frac{\lVert A\rVert}{\kappa})(\frac{2}{\frac{1}{l_{n+1}}-\frac{1}{l_{n+2}}} \lvert \ln \epsilon_{l_{n+1}} \rvert)^{\tau}\leqslant \epsilon_{l_{n+1}}^{-\frac{\sigma}{2}}.
$$

For the remaining estimates, we distinguish two cases.

If the $(n+1)$-th step is in the resonant case, we have
$$
A_{l_{n+1}}=e^{A_{l_{n+1}}''}, \ \ \lVert A_{l_{n+1}}''\rVert < 2\epsilon_{l_{{n+1}}}^{\sigma}, \ \ \lVert A_{l_{n+1}}\rVert\leqslant 1+2\epsilon_{l_{{n+1}}}^{\sigma}\leqslant 2\lVert A\rVert.
$$
Then there exists unitary $U\in SL(2,\C)$ such that
\begin{equation}
\label{estt}UA_{l_{n+1}}U^{-1}=\begin{pmatrix} e^{\gamma_{n+1}} & c_{n+1}\\ 0 & e^{-\gamma_{n+1}} \end{pmatrix},
\end{equation}
with $\lvert c_{n+1}\rvert\leqslant 2\lVert A_{l_{n+1}}''\rVert \leqslant 4\epsilon_{l_{{n+1}}}^{\sigma}$. Thus $(\ref{estisharp})$ is fulfilled because
\begin{equation}\label{non1}
\lVert B_{l_{n+1}}(\theta) \rVert_0^2\lvert c_{n+1}\rvert\leqslant 4.
\end{equation}

If it is in the non-resonant case, one traces back to the resonant step $j$ which is closest to $n+1$.

If $j$ exists, by $(\ref{estima})$ and $(\ref{estres})$ we have
$$
\lvert B_{l_j}(\theta)\rvert_{\frac{1}{l_{j+1}}}\leqslant \epsilon_{l_j}^{-\frac{\sigma}{2}-s},\ \ \lVert B_{l_j}(\theta)\rVert_0\leqslant \epsilon_{l_j}^{-\frac{\sigma}{2}},
$$
$$
A_{l_j}=e^{A_{l_j}''}, \ \ \lVert A_{l_j}''\rVert < 2\epsilon_{l_{j}}^{\sigma}, \ \ \lVert A_{l_j}\rVert \leqslant 1+2\epsilon_{l_{j}}^{\sigma}.
$$
By our choice of $j$, from $j$ to $n+1$, every step is non-resonant. Thus by $(\ref{estnon})$ we obtain
\begin{equation}
\label{estnn}\lVert A_{l_{n+1}}- A_{l_j}\rVert\leqslant 2\epsilon_{l_j}^{\frac{1}{2}},
\end{equation}
so
$$
\lVert A_{l_{n+1}}\rVert \leqslant 1+2\epsilon_{l_{j}}^{\sigma}+ 2\epsilon_{l_j}^{\frac{1}{2}}\leqslant 2\lVert A\rVert.
$$
Estimate $(\ref{estnn})$ implies that if we rewrite $A_{l_{n+1}}=e^{A_{l_{n+1}}''}$, then
$$
\lVert A_{l_{n+1}}''\rVert \leqslant 4\epsilon_{l_{j}}^{\sigma}.
$$
Moreover, by $(\ref{estb})$, we have
$$
\lVert B_{l_{n+1}}(\theta)\rVert_0 \leqslant \sqrt{2}\lvert B_{l_j}(\theta)\rvert_0 \leqslant \sqrt{2}\epsilon_{l_j}^{-\frac{\sigma}{2}}.
$$
Similarly to the process of $(\ref{estt})$, $(\ref{estisharp})$ is fulfilled because
\begin{equation}\label{non2}
\lVert B_{l_{n+1}}(\theta)\rVert_0^2 \lvert c_{n+1}\rvert\leqslant 8.
\end{equation}

If $j$ vanishes, it immediately implies that from $1$ to $n+1$, each step is non-resonant. In this case, $\lVert A_{l_{n+1}}\rVert \leqslant 2\lVert A\rVert$ and the estimate $(\ref{estisharp})$ is naturally satisfied as
$$\lVert B_{l_{n+1}}(\theta)\rVert_{\frac{1}{l_{n+1}}}\leqslant 2.$$
\end{pf}

Now, we are ready to show the quantitative almost reducibility theorem for $C^{k}$ quasi-periodic $SL(2,\R)$ cocycles.

\begin{Theorem}\label{thm3}
Let $\alpha\in {\rm DC}(\kappa,\tau)$, $\sigma<\frac{1}{6}$, $A\in SL(2,\R)$, $f\in C^k(\T^{d},sl(2,\R))$ with $k>(D+2)\tau+2$, there exists $\epsilon_1=\epsilon_1(\kappa,\tau,d,k,\lVert A\rVert,\sigma)$ such that if $\lVert f\rVert_k\leqslant \epsilon_1$ then $(\alpha , Ae^{f(\theta)})$ is $C^{k,k_0}$ almost reducible with $k_0\in\N$, $k_0\leqslant \frac{k-2\tau-1.5}{1+s}$. Moreover, if we further assume
$(\alpha , Ae^{f(\theta)})$ is not uniformly hyperbolic, then there exists $B_{l_j}\in C^{\omega}_{\frac{1}{l_{j+1}}}(2\T^{d}$,
$SL(2,\C))$, $A_{l_j}\in SL(2,\C)$, $\tilde{F}^{'}_{l_j}\in C^{k}(\T^{d}$,
$gl(2,\C))$, such that
$$
B_{l_j}(\theta+\alpha)(Ae^{f(\theta)})B^{-1}_{l_j}(\theta)=A_{l_j}+\tilde{F}^{'}_{l_j}(\theta)
$$
with $$\lVert B_{l_j}(\theta)\rVert_0\leqslant \epsilon_{l_j}^{-\frac{\sigma}{2}},\ \ \lVert \tilde{F}^{'}_{l_j}(\theta)\rVert_0\leqslant \epsilon_{l_j}^{\frac{1}{4}}$$
and
$A_{l_j}=\begin{pmatrix} e^{\gamma_j} & c_j\\ 0 & e^{-\gamma_j} \end{pmatrix}$ with estimate

\begin{equation}\label{ess}
\lVert B_{l_j}(\theta)\rVert_0^2 \lvert c_j\rvert\leqslant 8\lVert A\rVert,
\end{equation}
 where $\gamma_j\in i\R$ and $c_j \in \C$. \end{Theorem}

\begin{pf}
We first deal with the $C^0$ estimate. By Proposition $\ref{pro33}$, we have for any $l_j$, $j\in \N^+$:
$$
B_{l_j}(\theta+\alpha)(Ae^{f_{l_j}(\theta)})B^{-1}_{l_j}(\theta)=A_{l_j}e^{f_{l_j}^{'}(\theta)},
$$
thus
$$
B_{l_j}(\theta+\alpha)(Ae^{f(\theta)})B^{-1}_{l_j}(\theta)=A_{l_j}e^{f_{l_j}^{'}(\theta)}+B_{l_j}(\theta+\alpha)(Ae^{f(\theta)}-Ae^{f_{l_j}(\theta)})B^{-1}_{l_j}(\theta).
$$
Denote
\begin{equation}\label{quat1}
A_{l_j}+\tilde{F}_{l_j}(\theta)=A_{l_j}e^{f_{l_j}^{'}(\theta)}+B_{l_j}(\theta+\alpha)(Ae^{f(\theta)}-Ae^{f_{l_j}(\theta)})B^{-1}_{l_j}(\theta).
\end{equation}
In $(\ref{2.1})$, by a simple summation we get
\begin{equation}\label{quat2}
\lVert f(\theta)-f_{l_j}(\theta)\rVert_0\leqslant \frac{c}{(2\lVert A\rVert)^{\tilde{D}}l_1^{D\tau+\frac{1}{2}}l_{j}^{k-1}},
\end{equation}
and
\begin{equation}\label{quat3}
\lVert f(\theta)\rVert_0+\lVert f_{l_j}(\theta)\rVert_0\leqslant \frac{c}{(2\lVert A\rVert)^{\tilde{D}}M^{D\tau+\frac{1}{2}}}+\frac{c}{(2\lVert A\rVert)^{\tilde{D}}C'M^{D\tau+\frac{1}{2}}}.
\end{equation}
Proposition $\ref{pro33}$ also gives the estimates
\begin{equation}\label{quat4}
\lVert B_{l_j}(\theta)\rVert_0\leqslant 64(\frac{\lVert A\rVert}{\kappa})(\frac{2}{\frac{1}{l_j}-\frac{1}{l_{j+1}}} \lvert \ln \epsilon_{l_j} \rvert)^{\tau} \leqslant \epsilon_{l_{j}}^{-\frac{\sigma}{2}}
\end{equation}
\begin{equation}\label{quat5}
\lvert f_{l_j}^{'}(\theta)\rvert_{\frac{1}{l_{j+1}}}\leqslant\epsilon_{l_j}^{3-\sigma},
\end{equation}
and
\begin{equation}\label{quat6}
\lVert A_{l_j}\rVert\leqslant 2\lVert A\rVert.
\end{equation}
Thus by $(\ref{quat1})-(\ref{quat5})$, we have
\begin{flalign}\label{Festi}
\lVert \tilde{F}_{l_j}(\theta)\rVert_0 \leqslant & \lVert A_{l_j}f_{l_j}^{'}(\theta)\rVert_0+\lVert B_{l_j}(\theta+\alpha)(Ae^{f(\theta)}-Ae^{f_{l_j}(\theta)})B^{-1}_{l_j}(\theta)\rVert_0 \\
\leqslant &\lVert A\rVert\epsilon_{l_j}^{3-\sigma}+(\frac{\lVert A\rVert^3}{\kappa^2})(\frac{2}{\frac{1}{l_j}-\frac{1}{l_{j+1}}} \lvert \ln \epsilon_{l_j} \rvert)^{2\tau}\times\frac{64\times64c}{(2\lVert A\rVert)^{\tilde{D}}l_1^{D\tau+\frac{1}{2}}l_{j}^{k-1}}\nonumber\\
\leqslant & \epsilon_{l_j}^{1+s} \nonumber.
\end{flalign}

Now let us prove estimate $(\ref{ess})$. Note that Proposition \ref{pro33} implies that we only need to rule out the possibility that $\gamma_j\in \R \backslash \{0\}$.

Assume that $spec(A_{l_j})=\{e^{\lambda_j},e^{-\lambda_j}\}, \lambda_j\in \R \backslash \{0\}$. If $\lvert\lambda_j\rvert>\epsilon_{l_j}^{\frac{1}{4}}$,
Assume that $spec(A_{l_j})=\{e^{\lambda_j},e^{-\lambda_j}\}, \lambda_j\in \R \backslash \{0\}$, then there exists $P\in SO(2,\R)$ such that
$$
PA_{l_j}P^{-1}=\begin{pmatrix} e^{\lambda_j} & c_j \\ 0 & e^{-\lambda_j}\end{pmatrix},
$$
with $\lvert c_j \rvert \leqslant \lVert A_{l_j}\rVert\leqslant 2\lVert A\rVert$.

If $\lvert\lambda_j\rvert>\epsilon_{l_j}^{\frac{1}{4}}$, Set $B=diag\{ \lVert 4A\rVert^{-3}\epsilon_{l_j}^{\frac{1}{4}}, \lVert 4A\rVert^{3} \epsilon_{l_j}^{-\frac{1}{4}} \}$, then
\begin{equation}\label{123}
BP(A_{l_j}+\tilde{F}_{l_j}(\theta))P^{-1}B^{-1}=\begin{pmatrix} e^{\lambda_j} & 0 \\ 0 & e^{-\lambda_j}\end{pmatrix}+ F(\theta),
\end{equation}
where $\lVert F(\theta)\rVert_0 \leqslant \frac{\epsilon_{l_j}^{\frac{1}{2}}}{C\lVert A\rVert^5}$.
We rewrite
$$
\begin{pmatrix} e^{\lambda_j} & 0 \\ 0 & e^{-\lambda_j}\end{pmatrix}+ F(\theta)=\begin{pmatrix} e^{\lambda_j} & 0 \\ 0 & e^{-\lambda_j}\end{pmatrix}e^{\tilde{f}(\theta)}
$$
with $\lVert \tilde{f}(\theta)\rVert_0\leqslant \frac{\epsilon_{l_j}^{\frac{1}{2}}}{C\lVert A\rVert^4}$.
Then by Remark $\ref{rem2}$ and Corollary 3.1 of \cite{houyou}, one can conjugate $(\ref{123})$ to
\begin{equation}\label{uh}
\begin{pmatrix} e^{\lambda_j} & 0 \\ 0 & e^{-\lambda_j}\end{pmatrix}\begin{pmatrix} e^{\tilde{f}^{re}(\theta)} & 0 \\ 0 & e^{-\tilde{f}^{re}(\theta)}\end{pmatrix}=\begin{pmatrix} e^{\lambda_j}e^{\tilde{f}^{re}(\theta)} & 0 \\ 0 & e^{-\lambda_j}e^{-\tilde{f}^{re}(\theta)}\end{pmatrix}
\end{equation}
with $\lVert \tilde{f}^{re}(\theta)\rVert_0 \leqslant \frac{\epsilon_{l_j}^{\frac{1}{2}}}{C\lVert A\rVert^2}$. Therefore $(\alpha , Ae^{f(\theta)})$ is uniformly hyperbolic, which contradicts our assumption. Now we only need to consider $\lvert\lambda_j\rvert \leqslant \epsilon_{l_j}^{\frac{1}{4}}$. In this case, we put $\lambda_j$ into the perturbation so that the new perturbation satisfies $\lVert \tilde{F}^{'}_{l_j} \rVert_0\leqslant \epsilon_{l_j}^{\frac{1}{4}}$ and
$$
A_{l_j}=\begin{pmatrix} 1 & c_j \\ 0 & 1\end{pmatrix}.
$$

Now let us deal with the differentiable almost reducibility. By Cauchy estimates, for $k_0 \in \N$ with $k_0\leqslant k $ and $n\in \N$ with $n\geqslant j$, we have
\begin{flalign*}
&\lVert B_{l_j}(\theta+\alpha)(Ae^{f_{{l_{n+1}}}(\theta)}-Ae^{f_{l_n}(\theta)})B^{-1}_{l_j}(\theta)\rVert_{k_0}\\
                            \leqslant &\sup_{\substack{
                             \lvert l\rvert \leqslant k_0,
                             \theta \in \T^{d}
                          }}\lVert (\partial_{\theta_1}^{l_1}\cdots\partial_{\theta_d}^{l_d})(B_{l_j}(\theta+\alpha)(Ae^{f_{l_{n+1}}(\theta)}-Ae^{f_{l_{n}}(\theta)})B^{-1}_{l_j}(\theta))\rVert \\
                          \leqslant & (k_0)!(l_{n+1})^{k_0}\lvert B_{l_j}(\theta+\alpha)(Ae^{f_{l_{n+1}}(\theta)}-Ae^{f_{l_{n}}(\theta)})B^{-1}_{l_j}(\theta)\rvert_{\frac{1}{l_{n+1}}}\\
                          \leqslant & (k_0)!(l_{n})^{(1+s)k_0}(\frac{\lVert A\rVert^3}{\kappa^2})(\frac{2}{\frac{1}{l_n}-\frac{1}{l_{n+1}}} \lvert \ln \epsilon_{l_n} \rvert)^{2\tau}{\epsilon_{l_n}}^{-\frac{\frac{4}{l_{n+1}}}{\frac{1}{l_n}-\frac{1}{l_{n+1}}}}\frac{64\times 64c}{(2\lVert A\rVert)^{\tilde{D}}l_1^{D\tau+\frac{1}{2}}l_{n}^{k-1}} \\
                          \leqslant & \frac{C_1}{l_n^{k-(1+s)k_0-2\tau-2s(D\tau+\frac{1}{2})-1}}
\end{flalign*}
where $C_1$ is independent of $j$.

By a simple summation we get
$$
\lVert B_{l_j}(\theta+\alpha)(Ae^{f(\theta)}-Ae^{f_{l_j}(\theta)})B^{-1}_{l_j}(\theta)\rVert_{k_0}\leqslant \frac{2C_1}{l_j^{k-(1+s)k_0-2\tau-2s(D\tau+\frac{1}{2})-1}}.
$$
Similarly by Cauchy estimates, we have
\begin{flalign*}
\lVert f_{l_j}^{'}(\theta)\rVert_{k_0}& \leqslant (k_0)!(l_{j+1})^{k_0}\lvert f_{l_j}^{'}(\theta)\rvert_{\frac{1}{l_{j+1}}}\\
&\leqslant (k_0)!(l_{j})^{(1+s)k_0}\times \epsilon_{l_j}^{3-\sigma}\\
&\leqslant (k_0)!(l_{j})^{(1+s)k_0}\times (\frac{c}{(2\lVert A\rVert)^{\tilde{D}} {l_j}^{D\tau+\frac{1}{2}}})^{3-\sigma}\\
&\leqslant \frac{C_2}{l_j^{(D\tau+\frac{1}{2})(3-\sigma)-(1+s)k_0}}
\end{flalign*}
where $C_2$ is independent of $j$.

Thus we have
\begin{flalign}\label{festi}
\lVert \tilde{F}_{l_j}(\theta)\rVert_{k_0} &\leqslant  \lVert A_{l_j}f_{l_j}^{'}(\theta)\rVert_{k_0}+\lVert B_{l_j}(\theta+\alpha)(Ae^{f(\theta)}-Ae^{f_{l_j}(\theta)})B^{-1}_{l_j}(\theta)\rVert_{k_0}\\
&\leqslant \frac{C_3}{l_j^{(D\tau+\frac{1}{2})(3-\sigma)-(1+s)k_0}}+\frac{2C_1}{l_j^{k-(1+s)k_0-2\tau-2s(D\tau+\frac{1}{2})-1}}. \nonumber
\end{flalign}
Note that although we write $k>(D+2)\tau+2$ in the assumption, we simply choose $k=[(D+2)\tau+2]+1$ in the actual operation process. So the quantity of $k$ is totally determined by $D$. Here ``$[x]$'' stands for the integer part of $x$.

So if $k_0\leqslant \frac{k-2\tau-\frac{3}{2}}{1+s}$, then
$$
\lVert \tilde{F}_{l_j}(\theta)\rVert_{k_0}\leqslant \frac{C_4}{l_j^{\frac{1}{6}}},
$$
which immediately shows that
$$\lim_{j \to +\infty}\lVert \tilde{F}_{l_j}(\theta)\rVert_{k_0}=0.$$
It means precisely that $(\alpha , Ae^{f(\theta)})$ is $C^{k,k_0}$ almost reducible.
This finishes the proof of Theorem \ref{thm3}.
\end{pf}

\begin{Remark}
In view of Corollary 3.1 in \cite{CCYZ}, we have $L(\alpha, Ae^{f(\theta)})=0$.
\end{Remark}

\begin{Remark}\label{re3.5}
The norm of the conjugation map $B_{l_j}$ can be adjusted easily by the variation of the parameters. More precisely, if we assume $D>\frac{t}{\sigma}$ with $t\geqslant 2$, then by the proof of Claim \ref{cl2} we have $\epsilon_{l_{n_{j+1}}}<\epsilon^t_{l_{n_{j}}}, \forall j\in \Z^+$ and
$$
\lVert B_{l_j}(\theta)\rVert_0\leqslant \epsilon_{l_j}^{-\frac{\sigma}{2t}\sum_{j=0}^{\infty}\frac{1}{t^j}}.
$$
This is quite useful for spectral applications. In certain cases, we need to reduce $B_{l_j}$ at the cost of enlarging $D$, i.e. the initial regularity $k$ increases.
\end{Remark}

In order to obtain the $\frac{1}{2}$-H\"{o}lder continuity of the Lyapunov exponent, we have to assume $D>\frac{5}{2\sigma}$ so that $\lVert B_{l_j}(\theta)\rVert_0\leqslant \epsilon_{l_j}^{-\frac{\sigma}{3}}$. Recall that the restriction on $\sigma$ is ``$\sigma<\frac{1}{6}$''. In order to make the initial regularity $k$ small, we need to fix $\sigma$ sufficiently close to $\frac{1}{6}$ in the beginning.

By Theorem \ref{thm3} and the Proof of Theorem 1.1 in \cite{CCYZ}, we have

\begin{Theorem}\label{thm1}
Let $\alpha \in {\rm DC}(\kappa,\tau)$, if $(\alpha, A)$ is $C^{k^{'},k}$ almost reducible with $k^{'}>k> 17\tau+2$, then for any continuous map $B:\T^{d} \rightarrow SL(2,\C)$, we have
\begin{equation}\label{holder}
\lvert L(\alpha, A)-L(\alpha, B)\rvert \leqslant C\lVert B-A\rVert_0^{\frac{1}{2}},
\end{equation}
where $C$ is a constant depending on $d,\kappa,\tau,A,k$.
\end{Theorem}
\begin{pf}
The proof is almost the same as that in \cite{CCYZ}, with the only difference that the interval $I_j$ becomes: $C\epsilon_{l_j}^{\frac{1}{4}}\leqslant \epsilon \leqslant c{\epsilon_{l_j}^{\frac{2}{9}}}$. Here $C,c$ are two constants depending on $d,\kappa,\tau, A$. It is clear that all the small $\epsilon$ tending to zero can be covered by the interval $\{I_j\}_{j\geqslant 1}$ since $l_{j+1}=l_j^{1+s}$, $0<s\leqslant \frac{1}{6D\tau+3}$.
\end{pf}

Similarly, by Theorem \ref{thm3} and the Proof of Theorem 1.2 in \cite{CCYZ}, we have
\begin{Theorem}\label{cor1}
Let $\alpha \in {\rm DC}(\kappa,\tau)$, $V\in C^{k}(\T^{d},\R)$ with $k>17\tau+2$, then there exists $\lambda_{0}$ depending on $V,d,\kappa,\tau,k$ such that if $\lambda<\lambda_0$, then we have the following:
\begin{enumerate}
\item For any $E\in\R$, $(\alpha,S_E^{\lambda V})$ is $C^{k,k_0}$ almost reducible with $k_0\leq k-2\tau-2$.
\item  $N_{\lambda V,\alpha}$ is $1/2$-H\"older continuous:
$$
N(E+\epsilon)-N(E-\epsilon)\leqslant C_0\epsilon^{\frac{1}{2}}, \, \forall\, \epsilon>0, \,\forall \,E\in \R,
$$
where $C_0$ depends only on $d,\kappa,\tau,k$.
\end{enumerate}
\end{Theorem}
\begin{Remark}
In fact for $(1)$, by Theorem \ref{thm3} we only require $k$ to be larger than $14\tau+2$ since $D>\frac{2}{\sigma}$ is enough for almost reducibility, which gives Theorem \ref{thm1.2}.
\end{Remark}

\subsection{Applications to Schr\"{o}dinger cocycles}
For our purpose, we will apply our $C^k$ almost reducibility theorem on a special type of $C^k$ quasi-periodic $SL(2,\R)$ cocycles: Schr\"{o}dinger cocycle $(\alpha, S_E^{\lambda V})$, where
\begin{equation}
S_E^{\lambda V}(\theta)=\begin{pmatrix} E-\lambda V(\theta) & -1 \\ 1 & 0 \end{pmatrix}.
\end{equation}
For the sake of unification, let us rewrite $S_E^{\lambda V}(\theta)=Ae^{f(\theta)}$ where
$$
A=\begin{pmatrix} E & -1 \\ 1 & 0 \end{pmatrix}.
$$

In the following, the assumptions of Theorem \ref{thm3} are always fulfilled by assuming the small condition on $\lambda$ in Theorem \ref{main}. It gives that $(\alpha, S_E^{\lambda V})$ is $C^{k,k_0}$ almost reducible for all $E\in \R$, particularly for $E\in \Sigma$. Now let us divide $\Sigma$ into countable sets of energy $E$ in the following way. For $m\in \Z^+$, define
\begin{equation}\label{def1}
K_m=\{E\in\Sigma \mid (\alpha, S_E^{\lambda V}) \, \mbox{has a resonance at $m$-th step} \}.
\end{equation}
Moreover, define
$$
K_0=\{E\in\Sigma \mid (\alpha, S_E^{\lambda V}) \, \mbox{is reducible} \},
$$
then we have
$$
\Sigma=\displaystyle\bigcup_{m=0}^{\infty}K_m.
$$

When the resonance occurs, we can depict each $K_m$ more precisely by the rotation number of $(\alpha, S_E^{\lambda V})$ and we denoted it by $\rho(E)$ for convenience.

\begin{Lemma}\label{rotation}
Assume $E\in K_m,\, m\geqslant 1$, there exists $n\in \Z^d$ with $0<\lvert n\rvert \leqslant N_m$ such that
$$
\lvert 2\rho(E)-\la n,\alpha \ra\rvert_{\T} \leqslant 5\epsilon_{l_m}^{\sigma}.
$$
where $N_m=5l_{m}\ln\frac{1}{\epsilon_{l_m}}$.
\end{Lemma}
\begin{pf}
By Theorem \ref{prop1} and Theorem \ref{thm3}, if $E\in K_m$, then by the definition of the resonant case we have
\begin{equation}\label{11}
\lvert 2\rho(\alpha, A_{l_{m-1}})- \la n',\alpha \ra \rvert_{\T} \leqslant \epsilon_{l_m}^{\sigma},
\end{equation}
for some $n'\in \Z^d$ satisfying $0<\lvert n'\rvert\leqslant N=\frac{2}{\frac{1}{l_m}-\frac{1}{l_{m+1}}} \lvert \ln \epsilon_{l_m} \rvert$ (note that $A_{l_0}=A$). In addition, after doing this resonant step, by $(\ref{constant})$ we have
\begin{equation}\label{12}
\lvert \rho(\alpha, A_{l_{m}})\rvert \leqslant  \lVert A_{l_{m}}^{''}\rVert \leqslant 2\epsilon_{l_m}^{\sigma}.
\end{equation}

Moreover, formula $(\ref{degpro})$ gives
\begin{equation}\label{13}
\rho(E)+\frac{\la \deg B_{l_m},\alpha  \ra}{2}=\rho(\alpha, A_{l_m}+\tilde{F}_{l_m}(\theta)).
\end{equation}

By the properties of rotation number and $(\ref{Festi})$, there exists a numerical constant $c$ such that
\begin{equation}\label{14}
\lvert \rho(\alpha, A_{l_m}+\tilde{F}_{l_m}(\theta))-\rho(\alpha, A_{l_m})\rvert \leqslant c \lVert\tilde{F}_{l_m}(\theta)\rVert_0^{\frac{1}{2}}\leqslant c \epsilon_{l_m}^{\frac{1+s}{2}}.
\end{equation}

Denote $\rho(E)+\frac{\la \deg B_{l_m},\alpha \ra}{2}-\rho(\alpha, A_{l_m})=\ast$, then by $(\ref{11})$-$(\ref{14})$, we have
\begin{flalign*}
&\lvert \rho(E)+\frac{\la \deg B_{l_m},\alpha  \ra}{2}\rvert_{\T}\\
=&\lvert \rho(E)+\frac{\la \deg B_{l_m},\alpha \ra}{2}-\rho(\alpha, A_{l_m})+\rho(\alpha, A_{l_m})\rvert_{\T}\\
=&\lvert \ast+\rho(\alpha, A_{l_m})\rvert_{\T}\\
\leqslant & c \epsilon_{l_m}^{\frac{1+s}{2}}+2\epsilon_{l_m}^{\sigma}\\
\leqslant & \frac{5}{2}\epsilon_{l_m}^{\sigma}.
\end{flalign*}

By Claim \ref{cl2} and Remark \ref{re3.5}, for $t\geqslant 2$, we have
$$
\lvert \deg B_{l_m}\rvert\leqslant N\times\sum_{j=0}^{\infty}\frac{1}{t^j}\leqslant 5l_{m}\ln\frac{1}{\epsilon_{l_m}},
$$
if we denote $N_m=5l_{m}\ln\frac{1}{\epsilon_{l_m}}$, the result follows immediately.
\end{pf}

In order to make more preparations, we denote the transfer matrices by
\begin{equation}\label{transfer}
A_n(E,\theta)=\prod\limits_{j=n-1}^0 S_E^{\lambda V}(\theta+j\alpha).
\end{equation}
We will show that nice quantitative almost reducibility indicates nice control on the growth of $A_n$ on each $K_m, \, m\geqslant 1$ (as will be shown in Section 4, we do not need to estimate things on $K_0$ since it is transcendentally excluded in the proof of purely absolutely continuous spectrum).

\begin{Lemma}\label{chosen}
Assume $k>35\tau+2$, then for every $E\in K_m, m\geqslant 1$, we have $\sup_{0\leqslant s\leqslant c\epsilon_{l_m}^{-1+\frac{\sigma}{2}}}\lVert A_s\rVert_0\leqslant \epsilon_{l_m}^{-\frac{2\sigma}{9}}$, where $c$ is a universal constant.
\end{Lemma}

\begin{pf}
For $E\in K_m, m\geqslant 1$, $(\alpha, S_E^{\lambda V})$ has a resonance at $m$-th step. By the resonant estimates $(\ref{tildea})$-$(\ref{constm})$ of Theorem \ref{prop1}, Theorem \ref{thm3} and Remark \ref{re3.5}, if $D>\frac{11}{2\sigma}$ (thus we assume $k>35\tau+2$), then there exist $B_{l_{m}}\in C^{\omega}_{\frac{1}{l_{m+1}}}(2\T^{d},SL(2,\C))$, $A_{l_{m}}\in SL(2,\C)$ and $\tilde{F}^{'}_{l_{m}}\in C^{k_0}(\T^{d},sl(2,\C))$ such that
$$
B_{l_{m}}(\theta+\alpha)(Ae^{f(\theta)})B^{-1}_{l_{m}}(\theta)=A_{l_{m}}+\tilde{F}^{'}_{l_{m}}(\theta),
$$
with
$$
\lVert B_{l_{m}}(\theta)\rVert_0\leqslant \epsilon_{l_{m}}^{-\frac{\sigma}{9}}, \, A_{l_{m}}=\begin{pmatrix} e^{\gamma_{m}} & 0\\ 0 & e^{-\gamma_{m}} \end{pmatrix}, \lVert \tilde{F}^{'}_{l_{m}}(\theta)\rVert_{0} \leqslant \epsilon_{l_{m}}^{1-\frac{\sigma}{2}}+\epsilon_{l_{m}}^{1+s},
$$
where ``$\epsilon_{l_{m}}^{1-\frac{\sigma}{2}}$'' is the upper bound of the off-diagonal's norm, see $(\ref{D})$ and ``$\epsilon_{l_{m}}^{1+s}$'' corresponds to the quantity of the perturbation, see $(\ref{Festi})$. Moreover, we have $\gamma_{m}\in i\R$.

Then we easily conclude
$$
\sup_{0\leqslant s\leqslant c\epsilon_{l_m}^{-1+\frac{\sigma}{2}}}\lVert A_s\rVert_0\leqslant \lVert B_{l_{m}}(\theta)\rVert_0^2 \leqslant \epsilon_{l_m}^{-\frac{2\sigma}{9}}.
$$
\end{pf}

\section{Spectral application: absolutely continuous spectrum}
With the dynamical estimates in hand, we can prove our main theorem easily. Let us cite the well known result shown by Gilbert-Pearson \cite{GP}.

\begin{Theorem}\label{gp}\cite{GP}
Let $\mathcal{B}$ be the set of $E\in\R$ such that the cocycle $(\alpha,S_E^{V})$ is bounded. Then $\mu_{V,\alpha,\theta}\vert \mathcal{B}$ is absolutely continuous for all $\theta\in\R$.
\end{Theorem}

Besides, let us recall two convenient results proved by Avila \cite{A01}.
\begin{Theorem}\label{avi}\cite{A01}
We have $\mu(E-\epsilon, E+\epsilon)\leqslant C\epsilon\sup_{0\leqslant s\leqslant C\epsilon^{-1}}\lVert A_s\rVert_0^2$, where $C>0$ is a universal constant.
\end{Theorem}

\begin{Theorem}\label{avi2}\cite{A01}
If $E\in \Sigma$ then for $0<\epsilon<1$, $N(E+\epsilon)-N(E-\epsilon)\geqslant c\epsilon^{\frac{3}{2}}$, where $c>0$ is a universal constant.
\end{Theorem}

\begin{Remark}\label{rm4.1}
The proof of Theorem \ref{avi2} requires the $\frac{1}{2}$-H\"{o}lder continuity of the integrated density of states in our case, which is ensured by Theorem \ref{cor1}.
\end{Remark}

\subsection{Proof of Theorem \ref{main}} \begin{pf} Denote $\mathcal{B}$ the set of $E\in\Sigma$ such that $(\alpha,S_E^{\lambda V})$ is bounded. Denote $\mathcal{R}$ be the set of $E\in\Sigma$ such that $(\alpha,S_E^{\lambda V})$ is reducible. Then Theorem \ref{gp} ensure that we only need to prove $\mu(\Sigma\backslash \mathcal{B})=0$ for $\mu=\mu_{\lambda V,\alpha,\theta}$, $\theta\in \R$.

Note that $\mathcal{R}\backslash \mathcal{B}$ has only $E$ such that $(\alpha,S_E^{\lambda V})$ is reducible to parabolic. By the Gap Labeling Theorem \cite{BLT,MP}, for any $E\in \mathcal{R}\backslash \mathcal{B}$, there exists a unique $m\in \Z^d$ such that $2\rho(\alpha, S_E^{\lambda V})\equiv{\la m,\alpha\ra} \mod 2\pi\Z$, which shows $\mathcal{R}\backslash \mathcal{B}$ is countable. Moreover, if $E\in \mathcal{R}$, then any non-zero solution $H_{\lambda V,\alpha,\theta}u=Eu$ satisfies $\inf_{n\in \Z}\lvert u_n\rvert^2+\lvert u_{n+1}\rvert^2>0$. So there are no eigenvalues in $\mathcal{R}$ and $\mu(\mathcal{R}\backslash\mathcal{B})=0$. Therefore, it is enough to prove $\mu(\Sigma\backslash \mathcal{R})=0$.

Define $K_m, m\in \Z^+$ as in $(\ref{def1})$. Since $\Sigma\backslash \mathcal{R}\subset \lim \sup K_m$, we will show that $\sum\mu(\overline{K_m})<\infty$. Because by the famous Borel-Cantelli Lemma, $\sum\mu(\overline{K_m})<\infty$ implies $\mu(\Sigma\backslash \mathcal{R})=0$.

For every $E\in K_m$, let $J_m(E)$ be an open $\epsilon_m=C\epsilon_{l_m}^{\frac{2\sigma}{3}}$ neighborhood of $E$. By Lemma \ref{chosen},
$$
\sup_{0\leqslant s\leqslant C\epsilon_m^{-1}}\lVert A_s\rVert_0^2\leqslant  \epsilon_{l_m}^{-\frac{4\sigma}{9}}.
$$
Moreover, by Theorem \ref{avi}
$$
\mu(J_m(E))\leqslant C\epsilon_{l_m}^{-\frac{4\sigma}{9}}\lvert J_m(E)\rvert,
$$
where ``$\lvert \cdot \rvert$'' stands for the Lebesgue measure. Now we take a finite subcover $\overline{K_m}\subset \bigcup_{j=0}^r J_m(E_j)$. By refining this subcover, we can assume that any $x\in \R$ is contained in at most 2 different $J_m(E_j)$.

By Theorem \ref{avi2}, we have
$$
\lvert N(J_m(E))\rvert\geqslant c\lvert J_m(E)\rvert^{\frac{3}{2}}.
$$
By Lemma \ref{rotation}, if $E\in K_m$, then
$$
\lvert 2\pi N(E)-\la n,\alpha \ra\rvert_{\T} \leqslant 5\epsilon_{l_m}^{\sigma}
$$
for some $\lvert n\rvert\leqslant 5l_{m}\ln\frac{1}{\epsilon_{l_m}}$. This shows that $N(K_m)$ can be covered by $(10l_{m}\ln\frac{1}{\epsilon_{l_m}}+1)^d$ intervals $I_s$ of length $5\epsilon_{l_m}^{\sigma}$. Since $\lvert I_s\rvert\leqslant C\lvert N(J_m(E))\rvert$ for any $s$ and $E\in K_m$, there are at most $2C+4$ intervals $J_m(E_j)$ such that $N(J_m(E_j))$ intersects $I_s$. We conclude that there are at most $C(10l_{m}\ln\frac{1}{\epsilon_{l_m}}+1)^d$ intervals of $J_m(E_j)$. Then
$$
\mu(\overline{K_m})\leqslant \sum_{j=0}^{r}\mu(J_m(E_j))\leqslant C(10l_{m}\ln\frac{1}{\epsilon_{l_m}}+1)^d \times C\epsilon_{l_m}^{-\frac{4\sigma}{9}}\times C\epsilon_{l_m}^{\frac{2\sigma}{3}}\leqslant l_m^{-\frac{\tau}{5}},
$$
which gives
$$
\sum_{m}\mu(\overline{K_m})\leqslant C.
$$
This finishes the proof.
\end{pf}

As we have seen, the technical requirement of ``$k>35\tau+2$'' is really due to the necessity for ensuring the slower growth of the conjugation maps as well as the faster decay of the perturbations. For ``$k>14\tau+2$'', although we have the almost reducibility, there is no nice quantitative control of the conjugacy and the perturbation. Indeed, understanding the competition between these two terms is essentially the core of using dynamical estimates to derive spectral results.

Last but not least, it is worth mentioning that our absolutely continuous result is itself a nice theorem to apply. In a future joint work with Silvius Klein, we are going to discuss the coexistence of spectral types in non-analytic classes.

\section{Acknowledgements}
The author would like to thank Jiangong You and Qi Zhou for useful discussions at Chern Institute of Mathematics, and is grateful to Pedro Duarte for his persistent support at University of Lisbon as well as to Silvius Klein for his consistent support from PUC-Rio. This work is supported by PTDC/MAT-PUR/29126/2017, Nankai Zhide Fundation and NSFC grant (11671192).

\end{document}